\documentclass[11pt]{amsart}
 
\theoremstyle{plain}
\newtheorem{thm}{Theorem}[section]
\newtheorem{theorem}[thm]{Theorem}

\newtheorem{lemma}[thm]{Lemma}

\usepackage{graphicx}

\numberwithin{equation}{section}

\title{Creating walls to avoid unwanted points in root finding and optimization} 
\author{Tuyen Trung Truong}
\date{\today}
\address{University of Oslo, 0316 Oslo, Norway }\email{tuyentt@math.uio.no}
\begin{document}

\begin{abstract}
\noindent 

In root finding and optimization, there are many cases where there is a closed set $A$ one likes that the sequence constructed by one's favourite method will not converge to A (here, we do not assume extra properties on $A$ such as being convex or connected).  For example, if one wants to find roots, and one chooses initial points in the basin of attraction for 1 root  $z^*$ (a fact which one may not know before hand), then one will always end up in that  root. In this case, one would like to have a mechanism to avoid this  point $z^*$ in the next runs of one's algorithm. 

Assume that one already has a method IM for optimization (and root finding) for non-constrained optimization. We provide a simple modification IM1 of the method to treat the situation discussed in the previous paragraph. If the method IM has strong theoretical guarantees, then so is IM1. As applications, we prove two theoretical applications: one concerns finding roots of a meromorphic function in an open subset of a Riemann surface, and the other concerns finding local minima of a function in an open subset of a Euclidean space inside it the function has at most countably many critical points. 

Along the way, we compare with main existing relevant methods in the current literature. We provide several examples in various different settings to illustrate the usefulness of the new approach.  

\end{abstract}

\maketitle



\section{The problem, motivation and main result} Here we present the problem and motivation, a brief literature survey, two main theoretical results on finding roots of meromorphic functions in a given domain, and the plan for the remaining of the paper. 

\subsection{The problem and motivation}

Solving equations is an important task one usually encounters in research and applications. Some examples are: finding periodic points of a map, finding the trajectory of an object (e.g. robot) obeying a certain system of equations coming from physical laws. 

Ever since the time of Abel, Ruffini and Galois, it is clear that one cannot find precise roots of a simple polynomial in 1 complex variable, and hence must utilise approximative methods. In this paper, we will concentrate on iterative algorithms, which are very easy to implement and use in practice. 

One can treat solving a system of equations $F(x)=0$, where $F:\mathbf{R}^m\rightarrow \mathbf{R}^k$, as a global optimization problem by the following common trick. Define $f(x)=||F(x)||^2/2:~\mathbf{R}^m\rightarrow \mathbf{R}$. Then finding roots to $F=0$ is equivalent to the following two problems (which an iterative method may be able to solve simultaneously): 1) show that $\min _{x\in \mathbf{R}^m}f(x)=0$, and 2) find global minimizers for $f$. Even if $F$ has no root, finding global minimizers of $f$ still makes sense and has important applications (e.g. in the Least Square Fit problem in statistics). 

Since finding global minimizer is NP-hard, research on effective numerical methods for global optimization can be helpful for studying the problem P v.s. NP.  

This paper considers the following general problem: 

{\bf Problem:} Let $X$ be a (complete) metric space, and let $f:X\rightarrow \mathbf{R}$ a (smooth) cost function.  Assume that an iterative method IM, aiming to find global minimizers of $f$, is used, which has the property that whenever an initial point $x_0\in X$, then a sequence $x_{n+1}=IM(x_n,f)$ is generated, with $x_{n}\in X$ for all $n$. Can we have a way to construct sequences which avoid the set $A$ while still having a big chance to converge to global minima? In other words, can we use IM to solve the global optimization problem $f$ on the non-complete set $X\backslash{A}$?

For this problem to make sense, as a first requirement, the method IM must itself have a  big chance of convergence to global minimisers  on the complete set $X$. We note that well known methods like Gradient Descent and Newton's method, as well as many variants of them, do not have strong convergence guarantees, and hence need a lot of care when applied to functions which are non-convex and which have non-Lipschitz continuous gradients or Hessians. For example, a general polynomial of degree $4$ in 1 real variable does not satisfy those properties mentioned in the previous sentence.  

Below we list some examples for applying this problem. 

{\bf Application 1: Avoiding known roots.}  Assume that one wants to find roots to a system $F(x)=0$, with the associated cost function $f(x)=||F(x)||^2/2$. Assume that  in previous runs of IM one already found several roots $z_1^*,z_2^*,\ldots ,z_j^*$. One chooses the initial point for to run IM in a fixed domain $B$, and one does not know whether with probability 1, sequences constructed by IM with initial point in B will always converge to one of the known roots  $z_1^*,z_2^*,\ldots ,z_j^*$. Then one can define the close set $A=\{z_1^*,\ldots ,z_j^*\}$ and apply Problem. 

In the case where the variable x is in 1 dimension, an alternative way is to divide $F(x)$ by $(x-z_1^*)\times \ldots \times (x-z_j^*)$. However, one cannot do this in higher dimension. 

Even in dimension 1, there are situations when one would need to use Problem in its general setting. For example, one would like to find roots inside a certain domain $C$ (e.g. $C=\{z\in \mathbf{C}: ~|z|<1\}$). Then one can choose $A=\partial C$. 

{\bf Application 2: Avoiding known local minima/critical points.}  Assume that one wants to solve a global optimization problem for a cost function $f(x)$. Assume that  in previous runs of IM one already found local minima/critical points $z_1^*,z_2^*,\ldots ,z_j^*$. One does not know whether these are global minimizers (in particular, one does not want to stuck at bad local minima).  Then one can define the close set $A=\{z_1^*,\ldots ,z_j^*\}$ and apply Problem. 

{\bf Application 3: Constrained optimization.  } Assume that one has a function $f:X\rightarrow \mathbf{R}$, but only wants to solve the optimization problem in a smaller closed subset $Y$. Then one can choose $A$ to be the boundary of $Y$ and apply Problem.  

{\bf Application 4: Constrained optimization, version 2.} Assume that one has a function $f:X\rightarrow \mathbf{R}$, but only wants to solve the optimization problem in a smaller domain $Y$. Different from the choice in Application 3, here one chooses $A=X\backslash Y_0$, where $Y_0$ is the interior of $Y$.  

{\bf Application 5:  Finding different components of a variety.}  Assume the set of solutions to $F=0$ has different connected components $C_1,\ldots, C_k$, each of them  may have positive dimension. Assume that beforehand  one does not know any of them, and one will choose initial points for IM in a fixed domain $B$.  After many runs of IM, one finds many roots  $z_1^*,\ldots ,z_j^*$ but they seem to be close together and seem to belong to the same component $C_i$. (For example, this is the case if $B$ happens to belong to the basin of attraction of $C_i$). One hopes that if $j$ is big enough and one can avoid all of the points  $z_1^*,\ldots ,z_j^*$ , then one  can also avoid $C_i$ and finds another component of $\{F=0\}$.    

\subsection{Criteria for an optimization method to have strong theoretical guarantees} We clarify in this subsection our criteria for an optimization method IM to be regarded as having strong theoretical guarantees, which is an important input for our methods. Some explicit algorithms are discussed in some detail in the next section. 

Given an algorithm IM for optimization, we say that it has strong theoretical guarantees if it satisfies the following properties: (Here, $\{z_n\}$ is a sequence constructed by IM when applied to a function $f$ from an initial point $z_0$)

Property 1: Any cluster point of $\{z_n\}$ is a critical point of $f$. 

Property 2: There is a set $\mathcal{E}$ of Lebesgue measure $0$ so that if $z_0\notin \mathcal{E}$ and $\{z_n\}$ converges to $z^*$, then $z^*$ is not a saddle point of $f$. 

Property 3: For a large and useful class of functions $f$, the sequence $\{z_n\}$ will either converge to a finite point or diverge to $\infty$.  

Note that if Property 1--3 are satisfied,  then if $\{z_n\}$ has a bounded subsequence, then it must converge to a finite point $z^*$ which is a critical point of $f$ and which is not a saddle point of $f$. Hence, roughly speaking, $z^*$ should be a local minimum of $f$, as wanted in Optimization theory. This explains why we regard these properties as representing whether an optimization algorithm has strong theoretical guarantees or not.  

Here are a couple of examples of ``large and useful class of functions". The first class are those functions $f$ having at most countably many critical points (inside the domain of concern). This class includes the well known Morse functions, and is dense in the sense that for any function $f$ we can perturb it to become Morse. The second class are those functions $f$ satisfying the Lojasiewicz gradient inequality \cite{lojasiewicz}. This class includes  algebraic and analytic functions, as well as rational and meromorphic functions. These two classes cover almost all functions which one can encounter in practice.  

\subsection{A brief survey of main existing relevant methods in the current literature} The relevant methods which we are aware of are mostly in the setting of constrained optimization. 

{\bf Approach 1: New metric.} An approach is to redefine the metric on $X\backslash A$, such that the points in $A$ become infinity in the new metric. (This fits particularly well if $X$ is a Riemannian manifold.) However, if indeed the minimum of $f$ occurs in $A$, creating this new metric does not guarantee that the constructed sequence will not converge to $A$.  

{\bf Approach 2: Linear Programming.} If the cost function and all constraints are linear functions, then this is treated effectively in the well known Linear programming, see e.g. \cite{dantzig-thapa}. An experiment presented later will test whether our method can find good approximates of the optimizers for Linear Programming problems. 

{\bf Approach 3: Algebraic methods.} If the cost function and all constraints are polynomials, then there are algebraic methods aiming to solve it (e.g. \cite{lasserre}, \cite{pham}). We have a couple of remarks here. First, these methods are usually not iterative in nature (hence, implementation them in computers can be difficult). Second, they can - usually - only find the value $\min f(x)$ but not the points $z^*$ which minimize $f$. A combination between these and iterative methods may be more useful.

{\bf Approach 4: Projected methods.} If $X=\mathbf{R}^m$ and the constrained set $Y$ is a closed convex subset, then for each $x\in X$ there is a unique point $y\in Y$, denoted $y=pr_Y(x)$, so that $||x-y||=$ the distance from $x$ to $Y$. One can redefine the iterative scheme as follows: $x_{n+1}=pr_Y(IM(x_n,f))$.  However, if $A=\partial Y$ does not contain any global minimum of $f$, and if an open neighbourhood $U$ of $A$ belongs to the basin of attraction for a point in $X\backslash Y$, then choosing a random initial point  $x_0\in U$ and apply the projected method one ends up at a sequence in $A$ which can never converge to a global minimum of $f$. (An explicit example relating to this point was given in \cite{nouiehed-etal} . We will revisit this example in the experiments presented later.)   

{\bf Approach 5: Lagrange multiplier/Karush-Kuhn-Tucker conditions.} If the constraints are given by equations $h_1,\ldots ,h_j$, then Lagrange's multiplier method looks for critical points of a new function $F(x,\lambda _1,\ldots ,\lambda _j)=$ $f(x)-\lambda _1h_1(x)-\ldots \lambda _jh_j(x)$. This method works well if most of the critical points of $F$ corresponds to global minima of $f$, otherwise a lot of the work is wasted.  

Karush-Kuhn-Tuckter conditions \cite{karush}\cite{kuhn-tucker}: this extends Lagrange's method to the case where constraints also include inequalities. Then one looks for saddle points of a similar function. Again, the same comment as above for Lagrange's multiplier method can be applied here, and in this case  the issue may be more serious since saddle points are  more dominant in higher dimensions \cite{bray-dean}.  

{\bf Approach 6: Interior-point/Penalty methods.}  The main idea of these methods, see e.g. \cite{wright}, is to consider a new unconstrained problem $G(x)=f(x)-\epsilon \rho _A(x)$, where $\rho _A(x)$ is a function which is infinity on the boundary $A$, and $\epsilon >0$ is a parameter. It is unclear which $\epsilon$ is good to use. A usual choice of $\rho _A(x)$ is  $\rho _A(x)=\log d(x,A)$, where $d(.,.)$ is the distance function. The idea is that when $\epsilon$ becomes smaller and smaller, the method IM will find (for cost function $G(x)$) points which are closer and closer to global minimizers of $f(x)$.  

However, when $\epsilon $ is too small, what can happen is that it may cancel the effect of $ \rho _A(x)$  to a certain level that the constructed sequence by applying IM to $F(x)$ will behave similarly to the constructed sequence by applying IM to the original function $f(x)$. We will demonstrate this in an experiment later. 

{\bf Approach 7: Tunnelling/Deflation method.} This method deals with the case  where $A=\{z_1^*,\ldots ,z_j^*\}$ is a finite set. In the case of minimizing a function $f(x)$, it considers the new function $f(x)/[d(x,z_1)^{N_1}\ldots d(x,z_j)^{N_j}]$ for some appropriate  choices of $N_1,\ldots ,N_j$. It uses Backtracking line search to have descent property. In the case of solving a system of equations $F=0$, it applies Newton's method directly to the system, and again uses Backtracking line search to have descent property. Hence, it can show that the sequence constructed will avoid the finite set $A$. However, descent property alone does not guarantee strong convergence. For a constrained optimization problem, the tunnelling/deflation method does not work directly like our new methods below, but rather similar to the Lagrange multiplier/Karush-Kuhn-Tucker conditions. See \cite{levy-gomez}\cite{farrell-etal}\cite{farrell-etal2} for more details, where discussions on related methods are also given. The mentioned papers require the cost functions (or the system of equations) to have only discrete critical points (or roots), and do not prove global convergence of algorithms used. In an experiment later, we see that using this function $f(x)/[d(x,z_1)^{N_1}\ldots d(x,z_j)^{N_j}]$ may make it difficult to escape a component of the solution set to another component. 

\subsection{Application 1: finding roots of a meromorphic function in a given domain} Here we present a result which guarantees finding roots of a meromorphic function of 1 complex variable in a given domain. Similar results  in higher dimensions or for Riemann surfaces, and for functions satisfying the Lojasiewicz gradient inequality, can be stated with almost the same proof, by using the ideas in \cite{truong} and \cite{truong2}. 

Finding roots of a meromorphic function of 1 complex variable is a topic intensively studied. Many special, interesting and useful functions  are meromorphic functions, like Gamma function, Bessel function, and Riemann zeta function. The special question of finding roots of polynomial functions is a main subject in the fields of Complex Dynamics and Computer Science. 

Often, one has the need of finding roots in a given domain $U\subset \mathbf{C}$, for example inside the unit disk. Then, starting from a point in that domain $U$ and using an iterative method, it is difficult to know before hand if the sequence one constructed will converge to a root inside $U$. The reason is that, first of all the sequence constructed may not converge to any root at  all (in particular if one's method does not have good convergence guarantees), or it may converge to a root outside of $U$. Basins of attraction of the roots can be too complicated (e.g. that for Newton's method) to  allow a good guess of where the sequence will go. 

Here, we illustrate the use of the approach proposed in this paper towards this question. For the description and properties  of Backtracking New Q-Newton's method, see \cite{truong}. 

\begin{theorem} Let $g(z)$ be a non-constant meromorphic function in 1 complex variable $z$. Assume that $g$ is generic, in the sense that $\{z\in \mathbf{C}:~g(z)g"(z)=g'(z)=0\}=\emptyset$. 

Let $f(x,y)=|g(x+iy)|^2/2$. Let $U\subset \mathbf{C}$ be an open subset. There is a set $\mathcal{E}$ with Lebesgue measure 0 such that the following hold. 
 Let $M$ be a positive number and define a new function $F(x,y)$ by the following formula: $F(x,y)=f(x,y)$ if $(x,y)\in U$, and $F(x,y)=M$ if $(x,y)\in \mathbf{C}\backslash U$.  
 
 1) If $z_0=(x_0,y_0)\in U\backslash \mathcal{E}$, and $f(x_0,y_0)<M$, then Backtracking New Q-Newton's method applied to $F(x,y)$ with initial point $(x_0,y_0)$ will produce a sequence  $\{z_n=(x_n,y_n)\}$ which must satisfy one of the two options: a) $\{z_n\}$ converges to a root of $g(z)$ inside $U$ and the rate of convergence is quadratic, b)  all cluster points of $\{z_n\}$ are on $\partial U$. 
 
 2) Assume moreover that $U$ is bounded and $\inf _{z\in \partial U}f(z) >f(z_0)$. Then the constructed sequence (Backtracking New Q-Newton's method applied to F(x,y))  converges to a root of $g(z)$ inside $U$, with a quadratic rate of convergence. 
 
\label{Theorem}\end{theorem}
\begin{proof}
1) This follows easily from Theorem \ref{Theorem1} below, in combination with Theorem 3.3 in \cite{truong-etal}. 

2) Since $U$ is bounded, there is an open neighbourhood $V$ of $\partial U$ so that $\inf _{z\in V}f(z)>f(z_0)$. Therefore, given that the sequence $\{F(z_n)\}$ is non-increasing, no cluster point of $\{z_n\}$ can be on $\partial U$, hence option b) cannot happen. Hence, only option a) can happen. 
\end{proof}
Here are some remarks  on the condition $\inf _{z\in \partial U}f(z) >f(z_0)$ in the statement of the theorem.  Assume that $\inf _{z\in \partial U}f(z)>0$.  Then, by maximum principle, the necessary and sufficient condition for $g(z)$ to have a root inside $U$ is that there is a $z_0\in U$ for which $\inf _{z\in \partial U}f(z) >f(z_0)$. Therefore, our condition is optimal. 

\subsection{Application 2: Finding local minima of a function with countably many critical points in a domain} As another application, we consider the problem of finding local minima inside an open subset $U\subset \mathbf{R}^m$ for a function $f:\mathbf{R}^m\rightarrow \mathbf{R}$. We assume that $f$ has at most countably many critical points in $U$. We have an almost identical result like Theorem \ref{Theorem}. We can also extend it to other settings like Riemann manifolds, by using the ideas from \cite{truong, truong2}. 

 \begin{theorem} Let  $f:\mathbf{R}^m\rightarrow \mathbf{R}$ be a $C^3$ function. Let $U\subset \mathbf{R}^m$ be an open subset. Assume that $f$ has at most countably many critical points inside $U$. There is a set $\mathcal{E}$ with Lebesgue measure 0 such that the following hold. 
 Let $M$ be a positive number and define a new function $F(x)$ by the following formula: $F(x)=f(x)$ if $x\in U$, and $F(x)=M$ if $x\in \mathbf{R}^m\backslash U$.  
 
 1) If $z_0=\in U\backslash \mathcal{E}$, and $f(x_0)<M$, then Backtracking New Q-Newton's method applied to $F(x)$ with initial point $z_0$ will produce a sequence  $\{z_n\}$ which must satisfy one of the two options: a) $\{z_n\}$ converges to a critical point $z^*$ of $f(x)$ inside $U$, which is not a saddle point of $f$,  and if moreover $z^*$ is non-degenerate then the rate of convergence is quadratic, b)  all cluster points of $\{z_n\}$ are on $\partial U$. 
 
 2) Assume moreover that $U$ is bounded and $\inf _{z\in \partial U}f(z) >f(z_0)$. Then the constructed sequence (Backtracking New Q-Newton's method applied to F(x,y))  converges to a critical point of $f(z)$ inside $U$, which is not a saddle point of $f$, and if moreover $z^*$ is non-degenerate then the rate of convergence is quadratic. 
 
\label{TheoremB}\end{theorem}

\subsection{Acknowledgements and the plan of the paper} Several main ideas presented here were initiated in the inspiring environment of the ICIAM 2023 conference in Tokyo. The author would like to thank specially the organisers and participants of the minisymposium at ICIAM ``Theory and applications of random/non-autonomous dynamical systems" and its satellite conference, for  support and interesting discussions, in particular to Hiroki Sumi, Takayuki Watanabe and Mark Comerford. The author also would like to thank Xiao Wang on his informative presentation on convergence of iterative methods in a conference at University of Dalat in July 2023, and for pointing out the reference \cite{nouiehed-etal}.  We thank Patrick Farrell for information about the tunnelling/deflation method.  The author is partially supported by Young Research Talents grant 300814 from Research Council of Norway. 

The remaining of the paper is as follows. In Section 2 we present our two new methods. In Section 3 we present some illustrating experiments, including constrained optimization problems. In the last section we present some conclusions. There we also discuss some further ideas concerning constrained optimization. 

\section{Two new methods: Creating walls to avoid unwanted points} The main idea is to create walls at $A$, that is to make values of the functions at $A$ so large that it will unlikely for a sequence constructed from a point outside of $A$ can converge to or cross $A$. We propose here two new methods to create such walls. 

\subsection{Multiplying poles at $A$} The first method is to multiplying the cost function with a function having poles at $A$. 

We recall the setting: We are finding global minima of a function $f:X\rightarrow \mathbf{R}$, while wanting to avoid a closed subset $A\subset X$.  

We will present the method first in a special case, and then in the general case. 

\underline{\bf Special case: $\min _{x\in X}f(x)=0$}. 

In this case, we choose an appropriate positive integer $N$ and consider a new function $G(x)=f(x)/d(x,A)^N$, where $d(.,.)$ is the distance function. The number $N$ depends on the ``multiplicity'' of $f(x)$ along $A$. The following lemma  shows that at least we are on the right track with this method. 

\begin{lemma}
Assume that $f:X\rightarrow \mathbf{R}$ is such that $\min _{x\in X}f(x)=0$. Let $x^*\in X\backslash A$ be a root of $f(x)$ on $X$. Then $x^*$ is also a global minimizer of the new function $G(x)$. 
\label{Lemma1}\end{lemma}
\begin{proof}
Indeed, since $f(x)$ is non-negative, the same is true for $G(x)$. Since $x^*\in X\backslash A$ we have $d_A(x^*)\not= 0$. Hence, since $f(x^*)=0$, we have $G(x^*)=0$ as well. Therefore, $x^*$ is a global minimizer of $G$.  
\end{proof}

In the case of finding roots of a function $F(x)=0$ for $F:\mathbf{C}\rightarrow \mathbf{C}$, then this method is almost the same as what one usually does. Assume that one already found the roots $z_1^*,\ldots ,z_j^*$ of $F(x)$. Then the usual way to find new roots is to consider a new function $F(x)/[(x-z_1^*)\ldots (x-z_j^*)]$. In our method, we will try to find global minima of $f(x)=||F(x)||^2/2$, and will consider a  new cost function $G(x)=f(x)/d(x,A)^2$. Here $d(x,A) $ has the following explicit formula: $d(x,A)=\min \{|x-z_1^*|,\ldots ,|x-z_j^*|\}$.  Numerically, the complexity  of computing  $\nabla G$ and $\nabla ^2G$, as well as higher derivatives, is not increased when we increase the number of points $z_1,\ldots ,z_j$ to avoid.  

In the usual approach, if the function $f$ has multiplicity $N$ at a root $z_i^*$, one should divide by $(x-z_i^*)^N$. Similarly one should use $f(x)/d(x,A)^{2N}$ in our approach. [Note that since $f(x)=||F(x)||^2$, a root of multiplicity $N$ of $F$ will become a root of ``multiplicity" $2N$ for $f(x)$.] 

The tunnelling/deflation method is probably the one closest to our method of multiplying with poles at $A$. However, there are some differences. First, tunnelling/deflation method  only applies for the case where $A=\{z_1,\ldots , z_j\}$, while our method can apply to any closed set $A$. This is because we use a more flexible function $d(x,A)$, instead of the function $d(x,z_1)\ldots d(x,z_j)$. Second, while Backtracking line search is also used in tunnelling/deflation, it seems that they only use to guarantee descent property, and not Armijo's condition. Hence, while avoidance of the finite set $A$ is guaranteed, the convergence issue is not (except the case where some restrictive conditions are assumed). (Note that strong convergence guarantees for Armijo's Backtracking line search for general cost functions only appear quite recently in the literature.) Third, for to solve a system of equations $F=0$, the tunnelling/deflation method directly applies Newton's method for the system, while we work with the associated cost function $f=||F||^2$.    

Here we also note the relation between our method and the Penalty method recalled in the previous section. Since here the function $f(x)$ is non-negative, one can consider a new cost function $h(x)=\log f(x)$. Then, the Penalty method will require to consider the function $H(x)=\log f(x) - \epsilon \log d(x,A)$, where $\epsilon >0$ is {\bf small enough}. Taking logarithm of the function $G(x)=f(x)/d(x,A)^N$, we get another function $\log G(x)=\log f(x)-N \log d(x,A)$. Hence, the forms of $H(x)$ and $\log G(x)$ are the same, except a difference in detail: In $H(x)$ the factor $\epsilon $ before $\log d(x,A)$ is to be small, while in $\log G(x)$ the factor $N$ before $\log d(x,A)$ is to be big (at least $\geq 1$). Note also that generally our method using $G(x)$ is more numerically stable than using the Penalty method using $H(x)$, if IM at least uses gradients of the cost functions. Indeed, let $x^*\in X\backslash A$ be a point where $f(x^*)=0$  Then, near $x^*$ we have that $d(x,A)$ is strictly positive. Hence $\nabla G(x)$ is of the same size as $\nabla f(x)$, while $\nabla H(x)$ is of the same size as $(\nabla f(x))/f(x)$. Hence, near $x^*$, the gradient $\nabla H(x)$ is very big. Similar observation applies to points near $A$, where $d(x,A)$ is $0$.  An experiment showing the difference between these two approaches will be presented later. 

{\bf A heuristic argument.} While finding global minima is NP-hard, the following heuristic argument supports that our method will work well generally. 

Assume that $IM$ has the descent property for function values, that is whenever it is applied for a function $g$, then $g(x_{n+1})\leq g(x_n)$ for all $n$. There are several methods to ensure this property: Line search, Trust region, and Backtracking line search (Armijo or Frank-Wolfe).  Among these, Armijo's Backtracking line search is the most flexible while have strong theoretical guarantees when applied to many large and useful classes of cost functions. (For example, functions with at most countably many critical points,  or functions satisfying Lojasiewicz gradient inequality. The first class includes Morse functions, which is a dense class in the closed-open topology. The second class includes semi-analytic functions.) For some recent results on using Armijo's method on general functions, for both first and second order methods,  see \cite{truong-nguyen1}\cite{truong-nguyen2}\cite{truong-etal}\cite{truong}.  There, one also finds a brief comparison of several different methods.  

We also assume that with an appropriate choice of $N$, we obtain that $G(x)=\infty$ if and only if $x\in A$. 

Then for any choice of the initial point $x_0\in X\backslash A$, because IM has the descent property, it follows easily that the sequence $\{x_n\}$ constructed by IM from $x_0$ will belong to $X\backslash A$. Moreover, any cluster point of $\{x_n\}$ must be in $X\backslash A$. 

Now we explain why it is more likely that if $x_0$ is in a connected component of $X\backslash A$, then the whole sequence $\{x_n\}$ should be also in $X\backslash A$. We first present a rigorous argument, which is easy to prove. 

\begin{lemma}
Assume further that IM is such that there is a constant $r>0$ such that $d(x_n,x_{n+1})<r$ for all constructed sequence $\{x_n\}$. Let $A_r=\{x\in X:~d(x,A)\leq r\}$. Assume that $x_0\in X$ is such that $G(x_0)< \inf _{x\in A_r}G(x)$. Then the sequence $\{x_n\}$ constructed by IM will belong to the same connected component of $X\backslash A$ as $x_0$. 
\label{Lemma2}\end{lemma}
\begin{proof}
Indeed, assume by contradiction that some points in the sequence $\{x_n\}$ belongs to another component of $X\backslash A$. Let $N$ be the one index for which $x_N$ is still in the same connected component of $X\backslash A$ as $x_0$, but $x_{N+1}$ belongs to another component. Then since $d(x_N,x_{N+1})\leq r$, we must have $x_N\in A_r$. However, this contradicts the facts that IM is descent (hence, $G(x_N)\leq G(x_0)$) and $G(x_0)< \inf _{x\in A_r}G(x)$. 
\end{proof}

Even if the assumptions in  Lemma \ref{Lemma2} are not satisfied, we have the following heuristic argument. When we use Backtracking line search, then heuristically since on $A$ the function $G$ is $+\infty$, $A$ will play the roll of a wall which $\{x_n\}$ cannot cross. When $x_n$ is closer to $A$, the distance $d(x_n,x_{n+1})$ should become smaller in such a way that $x_{n+1}$ is always kept in the same connected component as $x_n$. 

This is how far heuristic arguments go, given that finding global minima are NP-hard. There are experiments which show that indeed there are cases where the sequence $\{x_n\}$ can cross to another connected component. Even so, we observe that by using the function $G(x)$, the chance for $\{x_n\}$ to stay in the same connected component indeed increases. 

\underline{\bf The general case}. We now consider a general function $f:X\rightarrow \mathbf{R}$, without assuming that it is non-negative (in particular, without assuming that $\min _{x}f(x)=0$). In theory, this  general case and the special case are equivalent, since if $f(x)$ achieves minimum on $X$, then we can replace $f$ by the function $f(x)-\min _{x}f(x)$ and one is reduced to the special case. However, finding $\min _{x}f(x)$ may not be easy, and hence in such cases one would like to be able to proceed without the hurdle to know that $\min _{x}f(x)=0$ (or even knowing if the function is non-negative). 

Some remarks are in order, before we continue. First, note that in this case the Penalty method which uses $\log f(x)$ is formally inapplicable, since the logarithm of a negative number is not defined (as a real number). Second, note that if the setting is semi-algebraic, one may use algebraic methods to find the minimum value of $f(x)$, and hence can reduce to the special case discussed before.

 To deal with this case where we do not know the minimum value of $f$, we consider the following function: $G_{\gamma}(x)=(f(x)-\gamma )/d(x,A)^N$, where $\gamma$ is an approximate lower bound for $min _{x}f(x)$ obtained ``on the fly". More precisely, we will proceed as follows. We will do many runs of the method IM for $G_{\gamma}(x)$, and in each run will try to obtain a better value for $\gamma$. In each run, we will do the following steps: 
 
 Step 1: Choose a random initial point $x_0$.
 
 Step 2: Apply IM to $G_{\gamma}(x)$ and construct a sequence $\{x_n\}$, where we will stop when some conditions are satisfied (for example, when $||\nabla G_{\gamma}(x_n)||$ is smaller than a prescribed threshold or when the number of iterates exceeds a given bound). 
 
 Step 3: If $\gamma > \min _n f(x_n)$, then we replace $\gamma$ by $ \min _n f(x_n)$. 
 
 If one wants to do optimization on a connected component $B$ of $X\backslash A$ only, then in Step 1 one chooses a random point $x_0\in B$, and one can modify Step 3 above in the following manner: 
 
 Step 3': If $\gamma > \min _{x_n\in B}f(x_n)$, then we replace $\gamma$ by $ \min _n f(x_n)$. 
 
 An example presented later will illustrate how this algorithm is used in practice. 
 
 \underline{\bf To escape the basin of attraction of a positive dimension component.}  We now consider the following question: Assume we apply  the iterative method IM to a non-negative cost function $f(x)$, whose zero set has different components $C_1,C_2,\ldots ,C_j$. If we found only points in $C_1$, can we escape it to reach another component? 
 
 Here we present an idea to use our method of multiplying poles, in combination of a careful choice of the next initial point, to resolve this question. 
 
 Assume that we already found a sequence of points (close to points) on $C_1$, called $p_1,\ldots ,p_k$.  As before, we will consider a new cost function $f(x)/d(x,A)^N$. Now, for the next initial point, we will not choose it randomly in a prescribed domain, but to close to the point $p_k$. This way, it can help us to quickly move away from $C_1$ and hopefully reach another component. 
 
 Here is a heuristic for why this idea can help. Because the component we attempt to escape has positive dimension, it can happen that even if we already created many poles, still in the next run we may end up at another point in that same component. When we creating more poles, we are also somehow creating a wall preventing initial points to escape to the other side. Hence, if the new poles are in between the domain where we choose our initial points and the component we want to reach, then may be we would not be able to reach that component as we wanted. To choose the initial point for the next pole close to the point (which now becomes a new pole) we found in the previous run, it is more likely that we avoid better the affect of the wall created by the poles.  
 
 An experiment later will illustrate how this is used. 
 
 \underline{\bf Constrained optimization.} This is a more precise form of the idea in Applications 3 and 4. In the above, we assumed that $G=+\infty$ on $A$ so that we can avoid $A$ (this is the case if $f>0$ on $A$). However, if $f<0$ at some points of $A$, the sequence we construct may converge to that point and we discover a possible minimum point. This way, we still can apply to constrained optimization. 
 
 If the constraints contain an equation $p(x)=0$, we replace it by two inequalties $p(x)\geq -\epsilon $ and $p(x)\leq \epsilon$, where $\epsilon >0$. The idea is that if $\epsilon$ goes to $0$, then the minimizers we found will converge to a minimizer of the original question. Numerically, with $\epsilon $ small enough, the minimizer we found for the new question should be very close to a minimizer of the original question.   
 
 Recall about the theorems for Backtracking New Q-Newton's method. 
 
 Note that if a function has compact sublevels, then when adding poles this way the new function will also has compact sublevels.  
 
 \subsection{Making the function to be a big constant on $A$} Another way to make walls at $A$ is to define the function $f$ to be a very large constant on $A$. More precise, we choose a very big number $R>0$, and then define a new function: $g(x)=f(x)$ if $x\in X\backslash A$, and $g(x)=R$ if $x\in A$. 
 
 While this method is very simple, it is guaranteed to avoid the interior of $A$. It also has the good property that it preserves critical points of $f$ inside $X\backslash A$ (this property, while simple, does not hold for all other approaches).  
 
 \begin{theorem}
 Assume that the iterative method IM has the descent property. Assume that IM is applied to $g$ at an initial point $x_0$ so that $f(x_0)<R$. Then no cluster point of the constructed sequence $\{x_n\}$ belongs to the interior of $A$. 
 
 The critical points (minima, maxima, saddle points) inside $X\backslash A$ of $f$ and $g$ are the same.    
 
 Any cluster point of $\{x_n\}$ which lies inside $X\backslash A$ is a critical point of $f$. 
 
 \label{Theorem1}\end{theorem}
 \begin{proof}
 Indeed, let $x_n$ be the constructed sequence. Then $f(x_n)\leq f(x_0)<R$ for all $n$. Hence $\{x_n\}\subset X\backslash A$, therefore no cluster point of $\{x_n\}$ belongs to the interior of $A$. 
 
 The last claim follows from that any cluster point of Backtracking New Q-Newton's method, applied to $f$, is a critical point of $f$.  
 \end{proof}
 Some experiments presented later illustrate that this method works well for constrained optimization.  Note that while the function $g$ is not continuous on the boundary $\partial A$ of $A$, usually $\partial A$ has zero Lebesgue measure and does not affect in numerical calculations.   
 
\underline {\bf Remark:} For constrained optimization, one might be tempted to think that the following approach - named Armijo's Backtracking line search with conditions - would produce the same effect as the above idea of making the function to be a big constant on $A$. 
 
 Armijo's Backtracking line search with constraints: Assume that one has a certain constraint $Constraint(x)$. In Armijo's Backtracking line search with constraints one uses the following while loop to find a good learning rate:
 
 {\bf While loop:} While (Armijo condition is not satisfied) or (Constraint(x) is not satisfied), replace learning rate $\kappa$ by  $\kappa /2$.
 
  However, an experiment later shows that the behaviour of these two approaches can be very much different. 
  
\subsection{Convergence issue} Since we will use Backtracking New Q-Newton's method, which has strong convergence guarantees and can avoid saddle points see \cite{truong-etal}\cite{truong}, as well as easy to implement and stable under perturbations, the issue of global convergence is less serious than if we use other methods, like Newton's method in the tunnelling/deflation method where a lot of care is needed. 
 
 \section{Experimental results} We now present several experimental results illustrating the usefulness of the approach proposed in this paper. 
 
 The setting of the experimental results is as follows. In a run for a cost function $g$, we will stop until either $||\nabla g||$ is smaller than a small threshold (e.g. $1e-6$) or if the number of iterates exceeds $10 000$. 
 
 In experiments where basins of attractions are drawn, we choose initial points on a square grid where the centre of the grid is chosen randomly in some domains. An initial point $x_0$ is considered to be in the basin of attraction of a point $z^*$ if $d(x_0,z^*)$ is smaller than a small threshold (e.g. $1e-5$). 
 
 Usually we  choose the initial point to be random, but there are cases where we fix the initial point. 
 
 Except Example 1- (which tests with methods Gradient Descent, Backtracking Gradient Descent (i.e. Armijo's Backtracking line search for Gradient Descent) and Backtracking New Q-Newton's method  \cite{truong} (developed from an earlier algorithm in \cite{truong-etal}) - all other examples test only Backtracking New Q-Newton's method. Backtracking New Q-Newton's method is found to be very much faster than the other two methods in the examples considered here. The exception is that Example 1 is taken from \cite{nouiehed-etal}, where Projected Gradient Descent was tested and found to fail to avoid the saddle point. Therefore, it is reasonable to test  Gradient Descent and Backtracking Gradient Descent for the function $G(x)$ in Example 1 also. 
 
 There is also a theoretical ground for our choice. Backtracking New Q-Newton's method, as mentioned above, has strong global convergence guarantee and can avoid saddle points \cite{truong}. More recent works \cite{fornaess-etal1, fornaess-etal2} show that Backtracking New Q-Newton's method has relation to the classical theorem of Schr\"oder in Complex Dynamics, as well as has connections to Newton's flow and Voronoi's diagrams. Gradient Descent, while can avoid saddle point, has no strong convergence guarantee. Projected Gradient Descent has less theoretical guarantees than Gradient Descent. Backtracking Gradient Descent, while has strong convergence guarantee, is not known whether can avoid saddle points. 
 
 \subsection{Example 1}  In this example, taken from \cite{nouiehed-etal}, we consider the following constrained optimization problem: 
 
$ arg\min _S f (x,y)$, where $S=\{(x,y)\in \mathbf{R}^2:~x+y\leq 0\}$ and $f(x,y)=-xye^{-x^2-y^2}+y^2/2$. 

One can check that the critical points of $f$ are (approximate to) the following: 

$p_1=(0,0)$, $p_2=(0.7071067,0.3128011)$ and $p_3=(-0.7071067,-0.3128011 )$. 

The point $(0,0)$ is a saddle point with function value $0$, while the other two points are global minimizer with function value (close to) $-0.0727279$. Therefore, $\gamma _0=-0.0727280$ is a good lower bound for $\min _Sf(x,y)$. 

In \cite{nouiehed-etal}, it was shown that Projected Gradient Descent applied to this constrained optimization problem fails. More precisely, there is a small open set $U$ in $S$ which touches the point $(0.5,-0.5)$ on the boundary for which if we apply Projected Gradient Descent  with a learning rate $0<\alpha <2/3$ with an initial point in $U$, then the constructed sequence will converge to the saddle point $(0,0)$. 

Here we draw basins of attraction for different methods for the function $G(x,y)=(f(x,y)-\gamma _0)/d((x,y),A)$ where $A=\partial S$. More explicitly, $d((x,y),A)=|x+y|$.  

We will also consider 2 functions $H_1(x,y)= (f(x,y)-\gamma _0) -\epsilon \log d((x,y),A)$, and $H_2(x,y)=\log  (f(x,y)-\gamma _0) -\epsilon \log d((x,y),A)$,  where $\epsilon >0$ is a constant.  $H_1$ is the Penalty method applied to the function $f(x,y)-\gamma _0$, and $H_2(x,y)$ is the Penalty method applied to the function $\log  (f(x,y)-\gamma _0) $. 
 
\begin{figure}
 \includegraphics[scale=1]{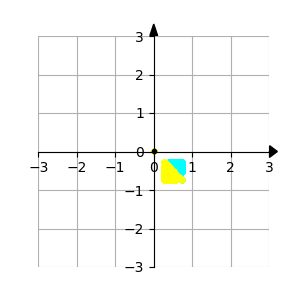}
\caption{Basins of attraction for the function G(x,y) in Example 1, using Gradient Descent with learning rate 0.1. Points are chosen on a square grid, with centre at the point $(0.5,-0.5003)$.  Cyan: initial points which  converge to $p_2$. Yellow: initial points which converge to $p_3$.} 
 \label{fig:pic1} 
\end{figure}

\begin{figure}
 \includegraphics[scale=1]{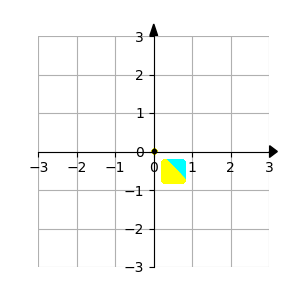}
\caption{Basins of attraction for the function G(x,y) in Example 1, using Backtracking Gradient Descent. Points are chosen on a square grid, with centre at the point $(0.5,-0.5003)$.  Cyan: initial points which  converge to $p_2$.  Yellow: initial points which converge to $p_3$. } 
 \label{fig:pic2} 
\end{figure}

\begin{figure}
 \includegraphics[scale=1]{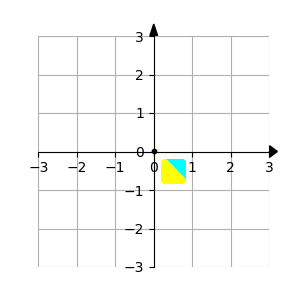}
\caption{Basins of attraction for the function G(x,y) in Example 1, using Backtracking New Q-Newton's method. Points are chosen on a square grid, with centre at the point $(0.5,-0.5003)$.  Cyan: initial points which  converge to $p_2$. Yellow: initial points which converge to $p_3$.} 
 \label{fig:pic3} 
\end{figure}

 \begin{figure}
 \includegraphics[scale=0.6]{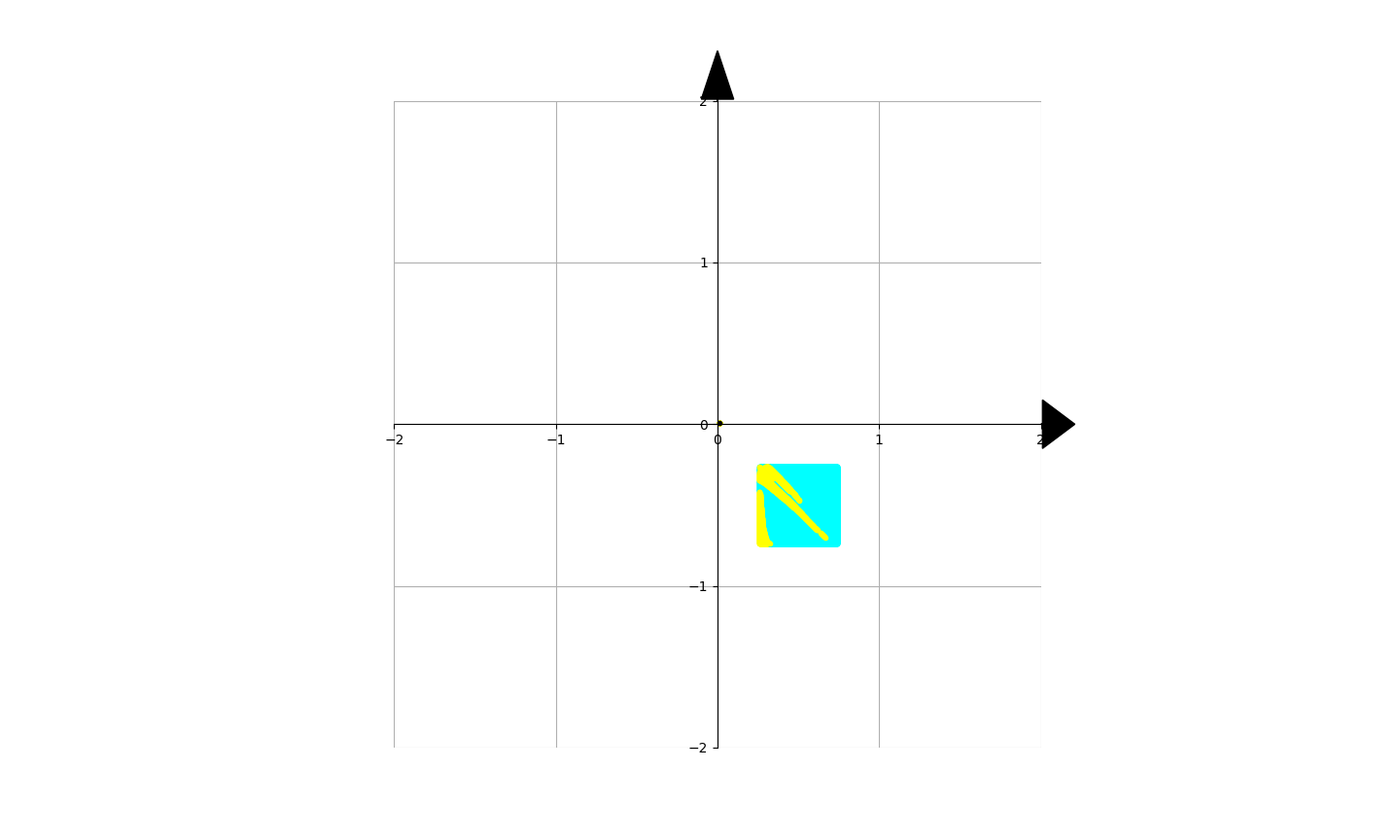}
\caption{Basins of attraction for the function $H_1(x,y)$ in Example 1, with $\epsilon =0.001$ using Backtracking New Q-Newton's method. Points are chosen on a square grid, with centre at the point $(0.5,-0.5003)$.  Cyan: initial points which  converge to $p_2$.  Yellow: initial points which converge to $p_3$. } 
 \label{fig:pic4} 
\end{figure}

 \begin{figure}
 \includegraphics[scale=0.6]{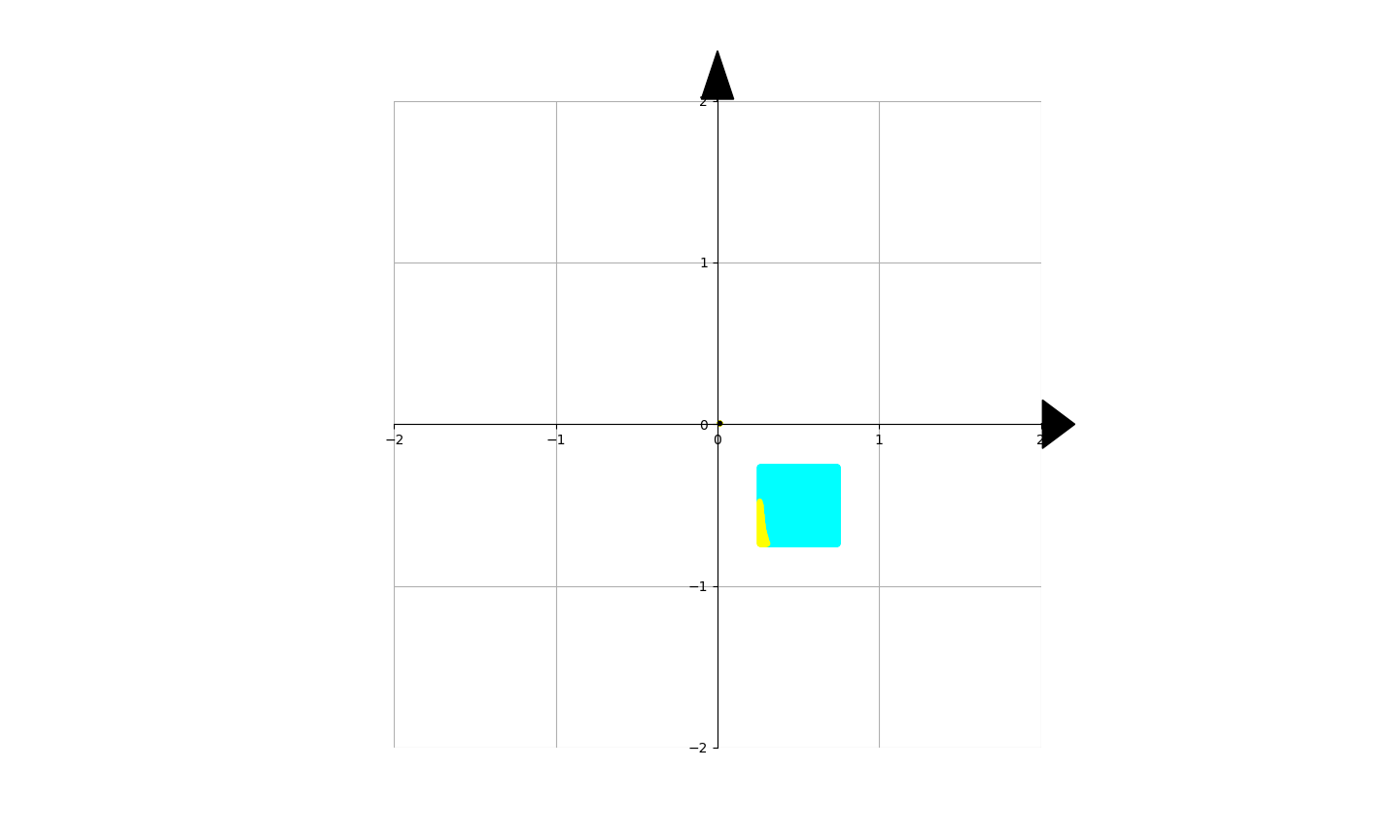}
\caption{Basins of attraction for the function $H_1(x,y)$ in Example 1, with $\epsilon =0.0001$ using Backtracking New Q-Newton's method. Points are chosen on a square grid, with centre at the point $(0.5,-0.5003)$.  Cyan: initial points which  converge to $p_2$.  Yellow: initial points which converge to $p_3$. } 
 \label{fig:pic4.1} 
\end{figure}

 \begin{figure}
 \includegraphics[scale=0.6]{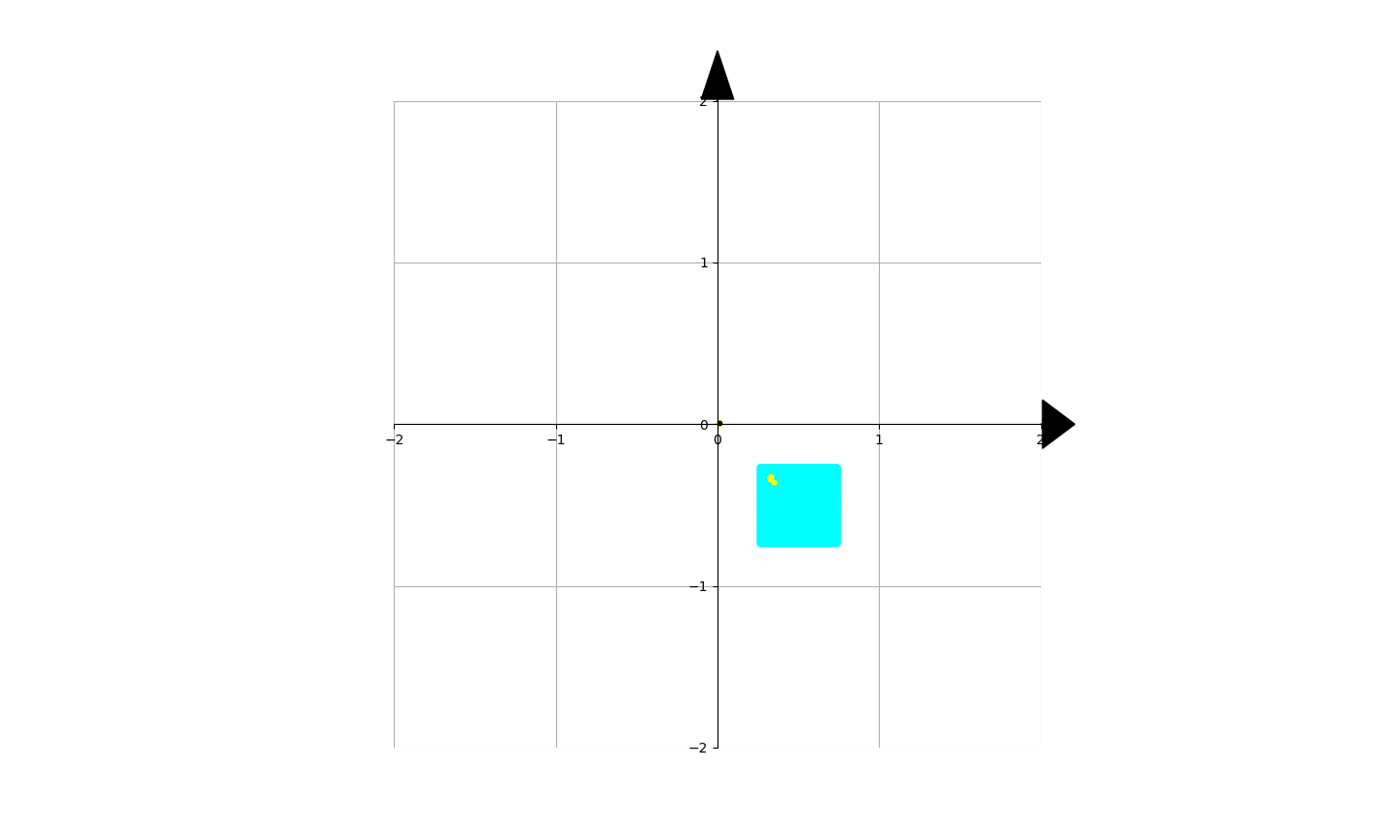}
\caption{Basins of attraction for the function $H_2(x,y)$ in Example 1, with $\epsilon =0.001$ using Backtracking New Q-Newton's method. Points are chosen on a square grid, with centre at the point $(0.5,-0.5003)$.  Cyan: initial points which  converge to $p_2$.  Yellow: initial points which converge to $p_3$. } 
 \label{fig:pic4.2} 
\end{figure}
 
 \begin{figure}
 \includegraphics[scale=0.6]{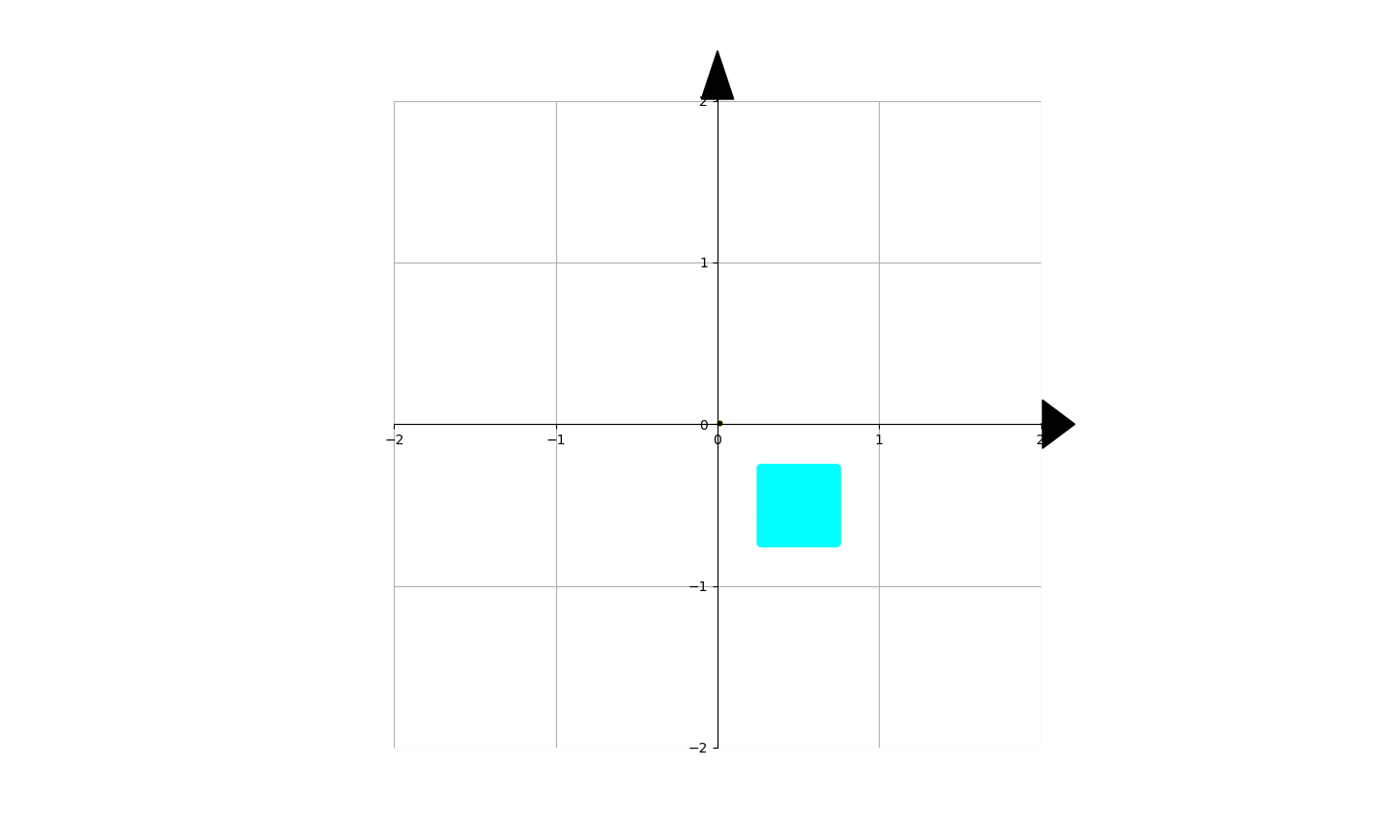}
\caption{Basins of attraction for the function $H_2(x,y)$ in Example 1, with $\epsilon =0.0001$ using Backtracking New Q-Newton's method. Points are chosen on a square grid, with centre at the point $(0.5,-0.5003)$.  Cyan: initial points which  converge to $p_2$.  Yellow: initial points which converge to $p_3$. } 
 \label{fig:pic4.3} 
\end{figure}

We have some remarks about experiments in Example 1. 

Figures \ref{fig:pic1}, \ref{fig:pic2} and \ref{fig:pic3} show that indeed using the function $G(x,y)$ helps improve the convergence to global minima. In particular, the phenomenon mentioned in \cite{nouiehed-etal} for Projected Gradient Descent does not seem to occur here.  The pictures for methods using Armijo's Backtracking line search look better than that for Gradient Descent.  

Figures \ref{fig:pic4},  \ref{fig:pic4.1},  \ref{fig:pic4.2} and  \ref{fig:pic4.3}  show that using Penalty method with $\log f(x,y)$, with small $\epsilon$ does not improve the performance. For $H_1(x,y)$, as we expected in the above, the performance will be like that for the function $f(x,y)$ when $\epsilon $ is small. Note that Figure  \ref{fig:pic4.1} is reminiscent of the phenomenon in Schroder's theorem for Newton's method for a polynomial of degree 2,  which is also numerically observed for Backtracking New Q-Newton's method see \cite{truong}, and proven rigorously in \cite{fornaess-etal1}. The yellow part in Figure  \ref{fig:pic4.1}  is close to the bisector of 2 points $p_2$ and $p_3$. Figures  \ref{fig:pic4.2} and  \ref{fig:pic4.3} show that, while formally optimising a non-negative function $f$ or its logarithm should be similar, there are enough differences (both theoretical and numerical), which lead to very different behaviours in experiments. We also remark that for $H_1$ and $H_2$, if we choose $\epsilon $ bigger, for example $\epsilon =0.01$, then we encounter errors such as NAN. 
 
\subsection{Example 2} In this example, we find roots of a polynomial of degree 5 in 1 complex variable. More precisely, the polynomial is $F(z)=z^5-3iz^3-(5+2i)z^2+3z+1$. This polynomial has the following (approximate roots): $p_1\sim -1.28992-1.87357i$, $p_2\sim -0.824853+1.17353i$, $p_3\sim -0.23744+0.0134729i$, $p_4\sim 0.573868-0.276869i$, and $p_5\sim 1.77834+0.963437i$. The associated cost function is $f(x,y)=|F(x+iy)|^2/2$, where $x,y\in \mathbf{R}$.  The following is a picture of basins of attractions when applying Backtracking New Q-Newton's method (it is drawn in a bigger domain than that in \cite{truong}).

 \begin{figure}
 \includegraphics[scale=0.6]{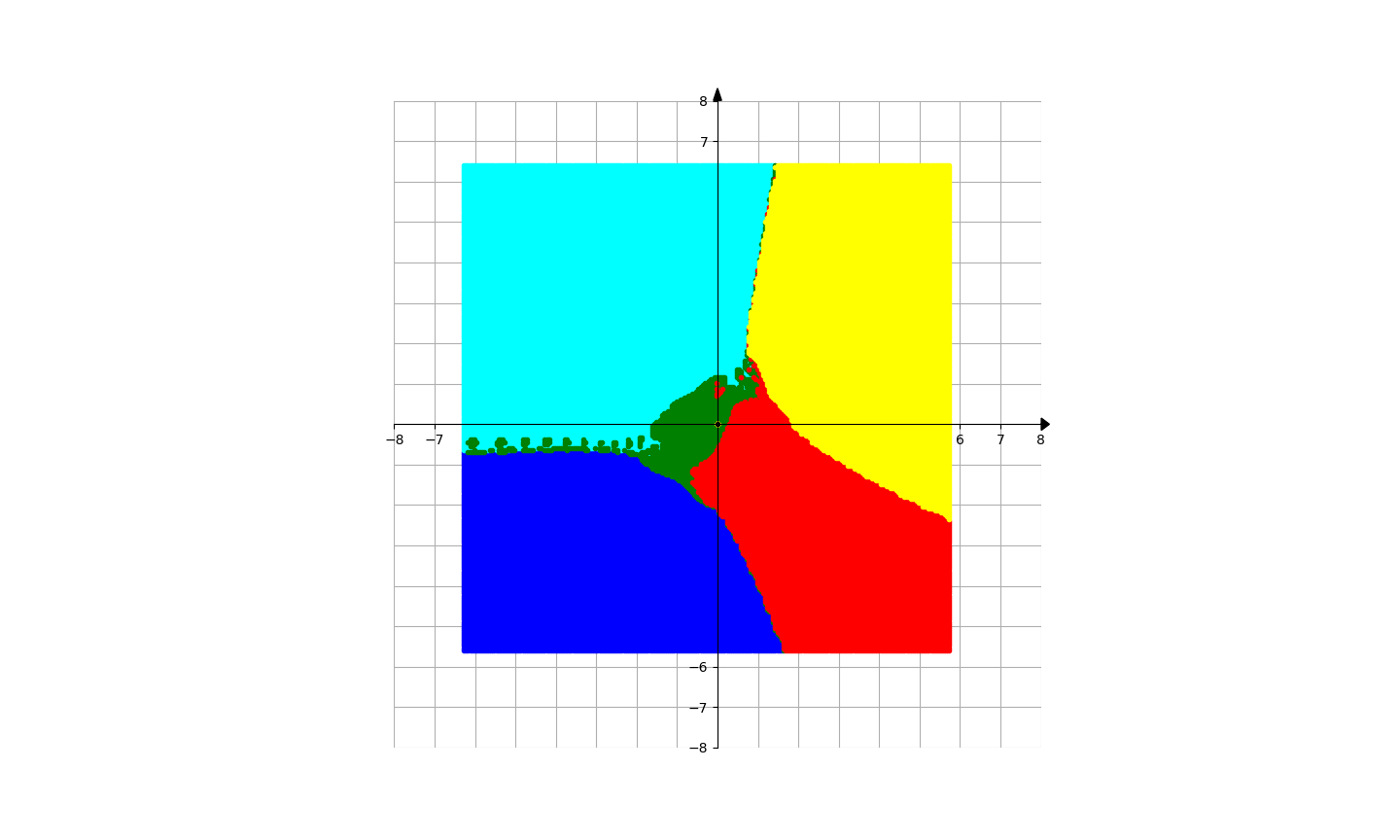}
\caption{Basins of attraction for the function f(x,y) in Example 2, using Backtracking New Q-Newton's method.  Blue: initial points which  converge to $p_1$.   Cyan: initial points which  converge to $p_2$.  Green: initial points which  converge to $p_3$. Red: initial points which  converge to $p_4$.   Yellow: initial points which converge to $p_5$. } 
 \label{fig:pic5} 
\end{figure}

We see from Figure \ref{fig:pic5} that the basin of attraction for $p_3$ seems to be smaller and surrounded by the basins of attraction of other points.  Hence, it may be difficult to find the root $p_3$. 

We will consider the following 4 closed sets: $A_1=\{p_1,p_2,p_4,p_5\}$, $A_2=\{p_1,p_2,p_4\}$, $A_3=\{p_1,p_2\}$ and $A_4=\{p_1\}$. Correspondingly, we consider 4 functions: $G_1(.)=f(.)/d(.,A_1)^2$,  $G_2(.)=f(.)/d(.,A_2)^2$, $G_3(.)=f(.)/d(.,A_3)^2$, and $G_4(.)=f(.)/d(.,A_4)^2$. 

 \begin{figure}
 \includegraphics[scale=0.3]{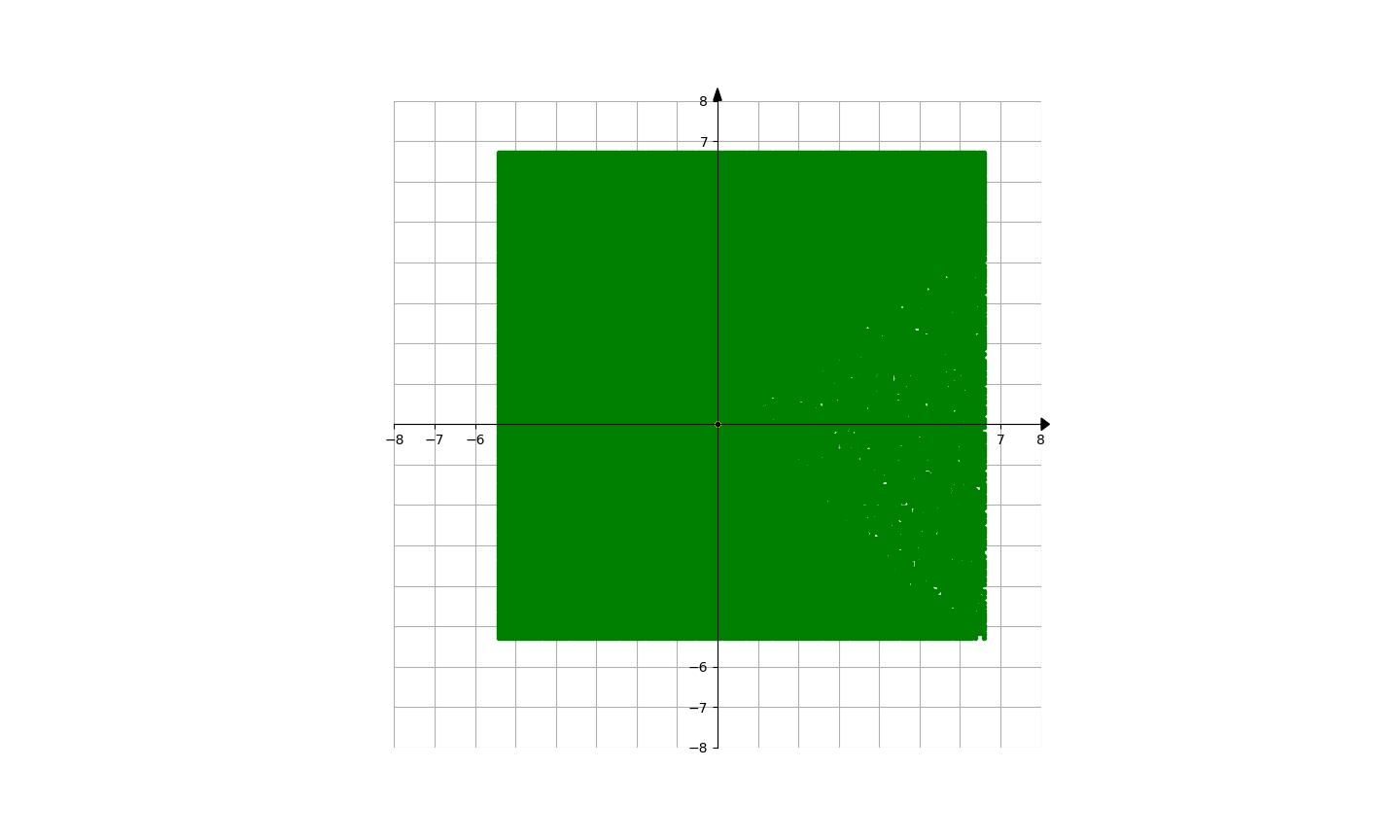}
\caption{Basins of attraction for the function $G_1(x,y)$ in Example 2, using Backtracking New Q-Newton's method.  Blue: initial points which  converge to $p_1$.   Cyan: initial points which  converge to $p_2$.  Green: initial points which  converge to $p_3$. Red: initial points which  converge to $p_4$.   Yellow: initial points which converge to $p_5$. } 
 \label{fig:pic6} 
\end{figure}

 \begin{figure}
 \includegraphics[scale=0.3]{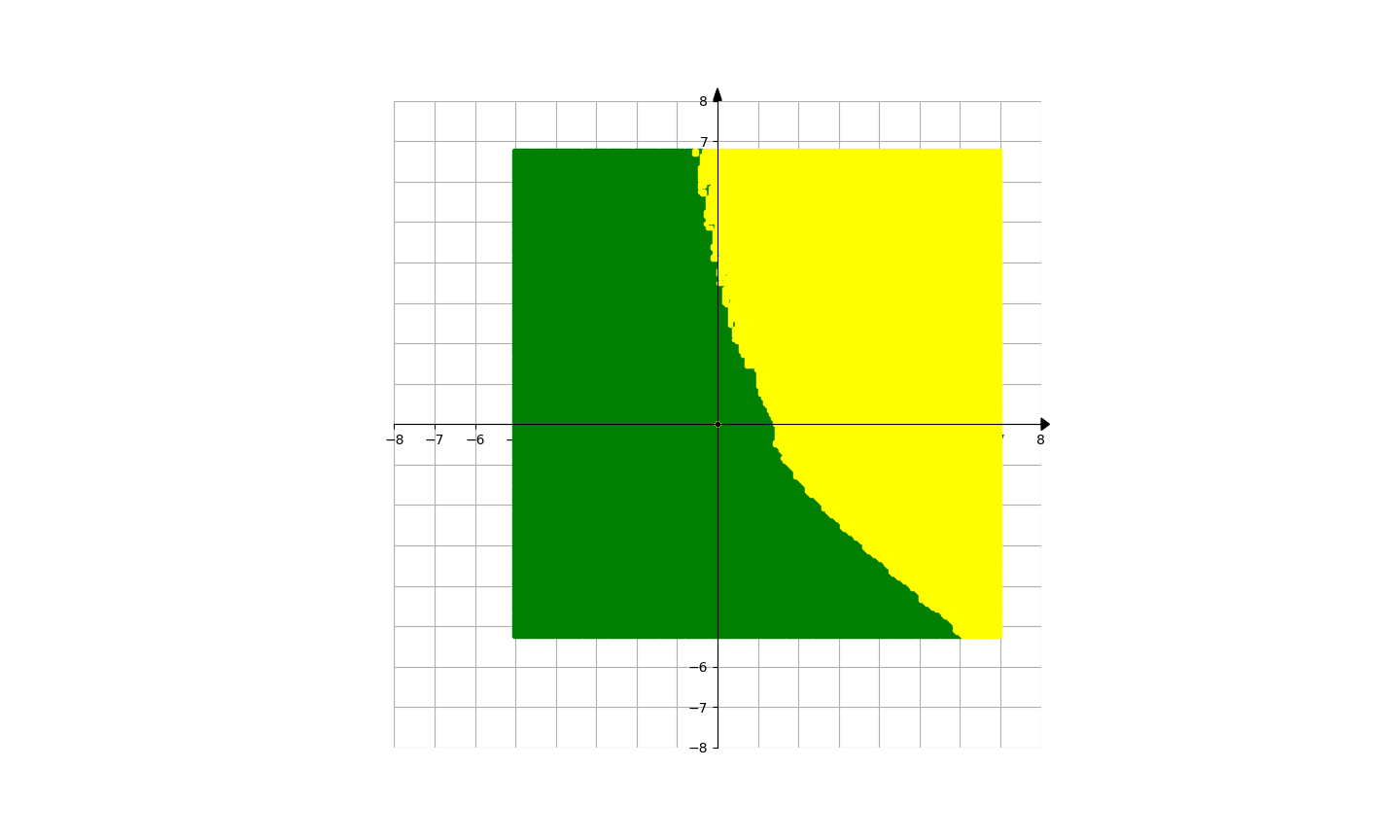}
\caption{Basins of attraction for the function $G_2(x,y)$ in Example 2, using Backtracking New Q-Newton's method.  Blue: initial points which  converge to $p_1$.   Cyan: initial points which  converge to $p_2$.  Green: initial points which  converge to $p_3$. Red: initial points which  converge to $p_4$.   Yellow: initial points which converge to $p_5$. } 
 \label{fig:pic7} 
\end{figure}

  \begin{figure}
 \includegraphics[scale=0.3]{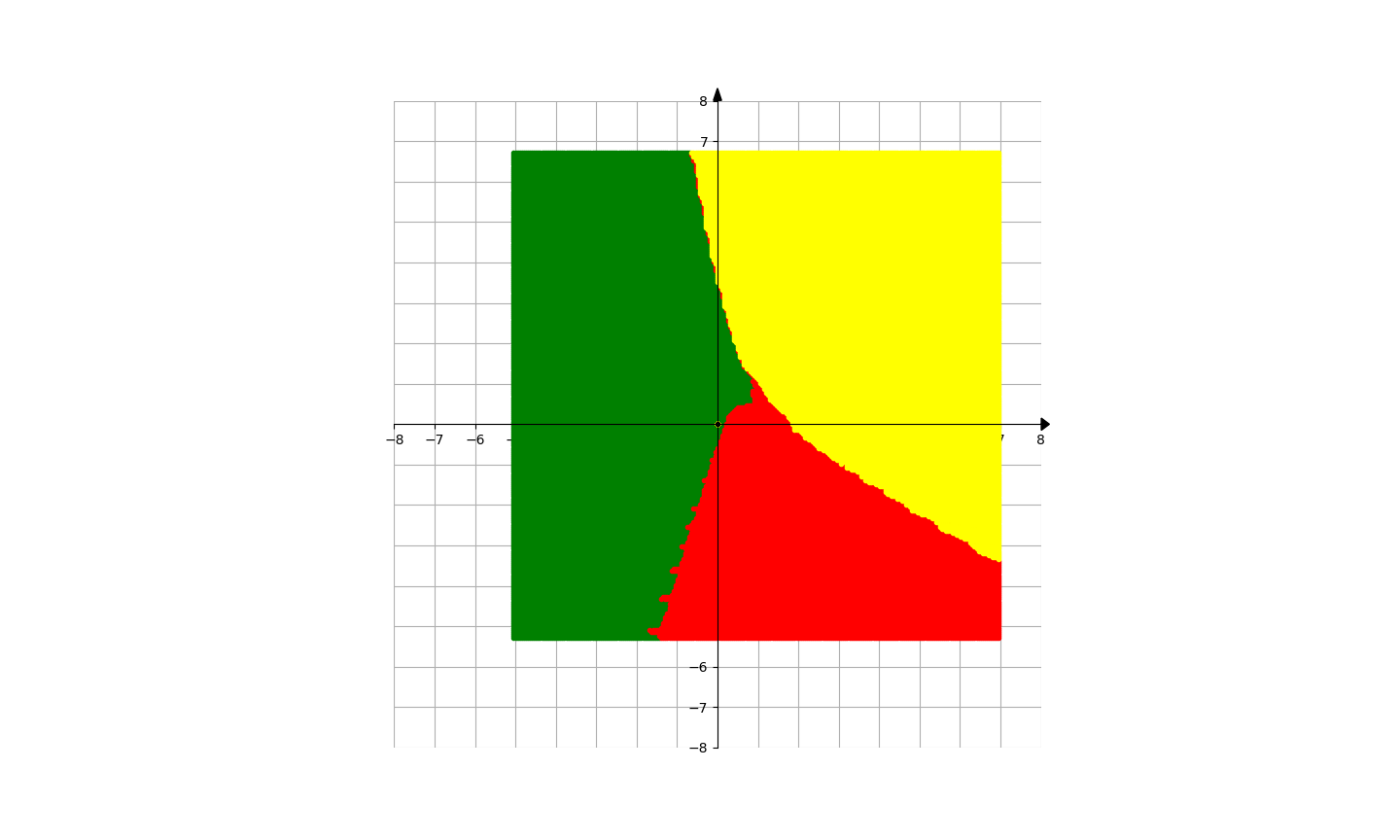}
\caption{Basins of attraction for the function $G_3(x,y)$ in Example 2, using Backtracking New Q-Newton's method.  Blue: initial points which  converge to $p_1$.   Cyan: initial points which  converge to $p_2$.  Green: initial points which  converge to $p_3$. Red: initial points which  converge to $p_4$.   Yellow: initial points which converge to $p_5$. } 
 \label{fig:pic8} 
\end{figure}

  \begin{figure}
 \includegraphics[scale=0.3]{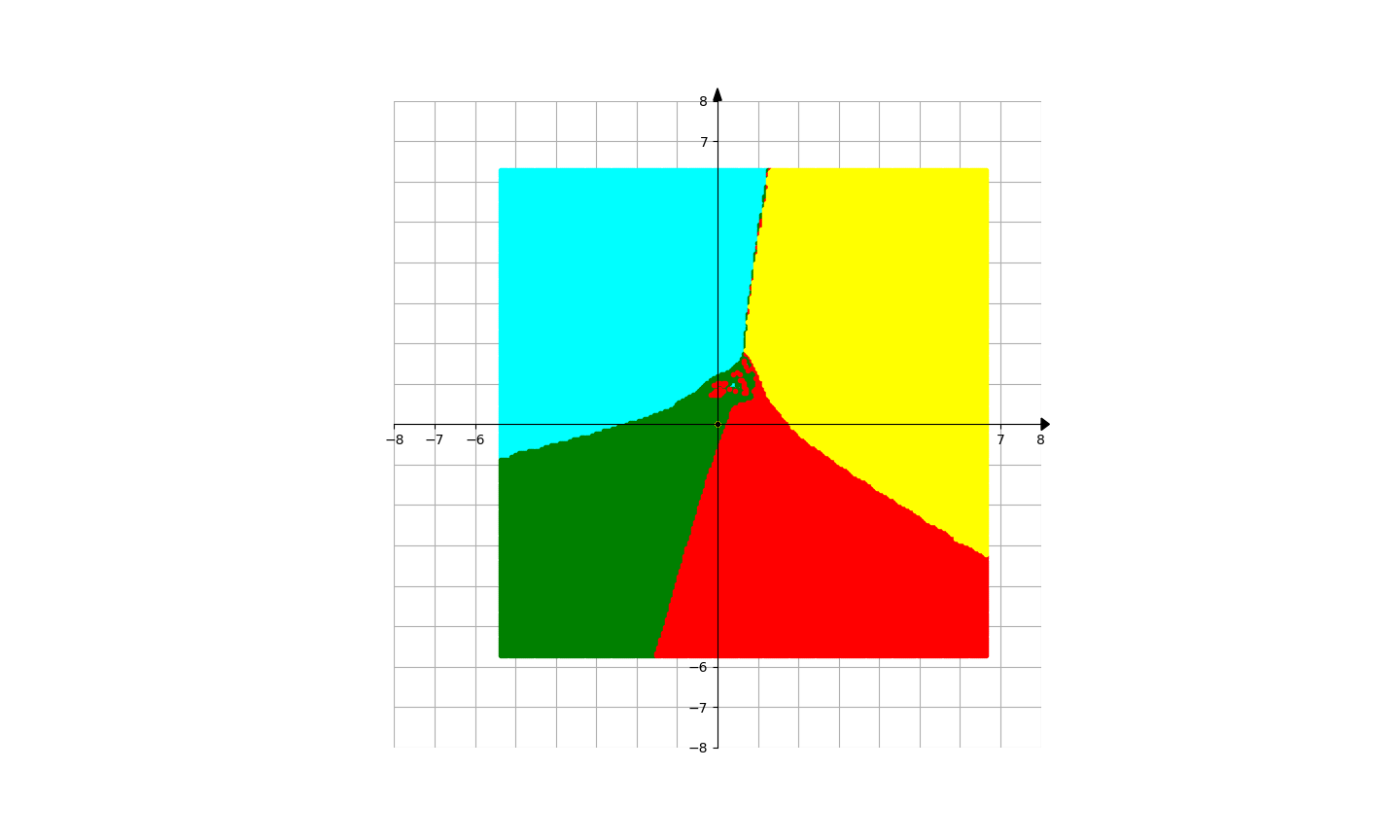}
\caption{Basins of attraction for the function $G_4(x,y)$ in Example 2, using Backtracking New Q-Newton's method.  Blue: initial points which  converge to $p_1$.   Cyan: initial points which  converge to $p_2$.  Green: initial points which  converge to $p_3$. Red: initial points which  converge to $p_4$.   Yellow: initial points which converge to $p_5$. } 
 \label{fig:pic9} 
\end{figure}

 \subsection{Example 3} We consider now finding roots of  a polynomial in 2 real variables $F(x,y)=(y-x^2-2)\times (y+x^4+2)\times (x^2+(y-1)^2)\times ((x-1)^2+(y+1)^2)\times ((x+1)^2+(y-4)^2)$. The solution set has 5 connected components, where 2 of them are curves: $C_1=\{(x,y)\in \mathbf{R}^2:~y-x^2-2=0\}$, $C_2=\{(x,y)\in \mathbf{R}^2:~y+x^4+2=0\}$, $C_3=(0,1)$, $C_4=(1,-1)$, $C_5=(-1,4)$. The first two components are curves, and the remaining components are points. Component $5$ is separated components $C_2,C_3$ and $C_4$ by component $C_1$. 
 
 In this experiment, we test if the initial point is in the set $y>x^2+2$, whether it is possible for an iterative use of the method in paper will allow to reach the components $C_2, C_3$ or $C_4$. 
 
 Since we want to solve equation, we choose our cost function as follows: $f(x,y)=(y-x^2-2)^2\times (y+x^4+2)^2\times (x^2+(y-1)^2)\times ((x-1)^2+(y+1)^2)\times ((x+1)^2+(y-4)^2)$
 
 a) We test with our method  using $f(x,y)/d((x,y),A)^2$. 
 
 We always choose the initial point to be $(0.1,3.1)$. 
 
 Run Backtracking New Q-Newton's method for $f(x,y)$, we end up at the point\\ $p_1=$ $(0.02299946, 2.00052898)$, which is close to a point on $C_1$. Set $A_1=\{p_1\}$. 
 
 Run Backtracking New Q-Newton's method for $f(x,y)/d((x,y),A_1)^2$, we end up at the point $p_2=(-0.2120813,2.04497848)$, which is again close to a point on $C_2$. Set $A_2=\{p_1,p_2\}$. 
 
 Run Backtracking New Q-Newton's method for $f(x,y)/d((x,y),A_2)^2$, we end up at the point $p_3=(-0.2559445,2.06550759)$,  which is again close to a point on $C_2$. Set  $A_3=\{p_1,p_2,p_3\}$.  
 
  Run Backtracking New Q-Newton's method for $f(x,y)/d((x,y),A_3)^2$, we end up at the point $C_3$.
  
  Thus, this example shows that it is possible to use our method iteratively to escape the basin of attraction for a big component and move to another component.   
  
b)  We next test if using $\prod _{p_j\in A}d((x,y),p_j)^2$, as in the tunnelling/deflation method, can help escape to another component. Again, we always choose the initial point to be $(0.1,3.1)$, which lies inside the domain $y>x^2+2$. 
  
   Run Backtracking New Q-Newton's method for $f(x,y)$, we end up at the point\\ $p_1=$ $(0.02299946, 2.00052898)$, which is close to a point on $C_1$. Set $A_1=\{p_1\}$. 

 Run Backtracking New Q-Newton's method for $f(x,y)/d((x,y),A_1)^2$, we end up at the point $p_2=(-0.2120813,2.04497848)$, which is again close to a point on $C_2$. Set $A_2=\{p_1,p_2\}$. 
 
  Run Backtracking New Q-Newton's method for $f(x,y)/\prod _{p_j\in A}d((x,y),p_j)^2$ $=$\\ $f(x,y)/[d((x,y),p_1)^2d((x,y),p_2)^2]$, we end up at the point $(-1,4)$, which is the component $C_5$. So we cannot escape the domain $y>x^2+2$.

   \subsection{Example 4} The test in this experiment is like that in Example 3, but now we look for roots of the real elliptic curve $E=\{y^2=x^3-x\}\subset \mathbf{R}^2$. This curve has two components $C_1=E\cap \{x\leq 0\}$ and $C_2=E\cap  \{x\geq 1\}$. Since we want to find roots, we choose the cost function to be $f(x,y)=(y^2-x^3+x)^2$. 
   
   We would like to start from a point close to $C_1$, and hope to reach $C_2$. 
   
   In this example, there are some differences to Example 3. 
   
   First, if we only divide by $d(.,A)^2$, it is not enough to escape $C_1$.  Hence, we change to $d(.,A)^4$.    
   
   Second, we find that if we choose the initial point in a prescribed small open set near $C_1$, we may not be able to reach $C_2$ but may wander around in $\mathbf{R}^2$. 
   
   Therefore, we will use the idea how to escape a component of positive dimension in  Section 0.3. 
   
   For example, assume that our random initial point happens to be $(-0.9, 0.1)$, which is inside the open set bound by $C_1$! After running Backtracking New Q-Newton's method for $f$, we end at the point $p=(-0.14223901, 0.3733089)$, which is (close to) a point on $C_1$.
   
   We now consider the new cost function $G(.)=f(.)/d(.,p)^4$. We now run Backtracking New Q-Newton's method for $G$, with the initial point $(-0.142, 0.373)$ which is close to the point $p$. We end up at the point $(1.01894051,-0.19739787)$, which is (close to) a point on the component $C_2$, as wanted!
   
\subsection{Example 5} In this experiment, we illustrate how the method in Section 0.3 applies for the case we do not know about $\min f(x)$ (and hence, not know if that value is actually $0$, which is the best case to use the method). 

We revisit the question in Example 1:  Consider the following constrained optimization problem: 
 
$ arg\min _S f (x,y)$, where $S=\{(x,y)\in \mathbf{R}^2:~x+y\leq 0\}$ and $f(x,y)=-xye^{-x^2-y^2}+y^2/2$.

Now, we do not before hand try to determine the critical points (and in particular, the minimum value) of the function $f$.  By the ideas in Section 0.3, we will start from $M=0$ and update $M$ to have better estimates for a lower bound for $\min f$ on the fly. 

Accordingly, the new cost function which we will apply Backtracking New Q-Newton's method to is $G_M(x,y)=(f(x,y)-M)/d((x,y),A)$, here $A=\{x+y=0\}\cup \{x=-10\} \cup \{x=10\} \cup \{y=1-0\} \cup \{y=10\}$. The reason to choose $A$ much bigger than the line $\{x+y=0\}$ is because it can easily check that the function $f$, while having only 3 critical points, has its gradient converges to $0$ for certain sequences of points going to $\infty$.  Hence (and this is observed in practice), the sequence constructed by an iterative algorithm can diverge to infinity.  To prevent this, we create new walls surrounding a bounded domain. 

Here, it turns out that always choosing a fixed point as the initial point for every run does not yield desired result. Hence, we switch to choosing our initial point randomly inside the domain $x+y<0$. Specifically, we will choose our random initial point in this way. We choose $x_0$ as a random point between $(-1,1)$ and also choose a random number $err\in (0,1)$. Then our initial point is $(x_0,y_0)=(x_0,-x_0-err)$.   

Beware that if $M\not= \min f(x)$, then global minimzer of $f$ may not be critical points of $G_M$! Therefore, applying an iterative method to $G_M$ may not be able to find a global minimizer of $f$. 

Here we report such an experiment. 

Start with $M_0=0$: The random initial point is $(0.23201923, -0.29366141)$. Running Backtracking New Q-Newton's method for $G_{M_0}(x,y)$, we get the sequence of points $(0.20221512,0.3141085)$, $(-0.04991277,  0.02411369)$, $(-0.05235699,  0.01419474)$, $(0.05467527, -0.01150241)$, \\ $(-0.70066302, -0.14982963)$, $(-0.52160544, -0.19009722)$, $(-0.51322439, -0.2418863)$, \\ $(-0.52258626, -0.24766466)$, $(-0.52258709, -0.24768905)$, $(-0.52258709, -0.24768905)$,\\ $(-0.52258709, -0.24768905)$, $(-0.5225871,  -0.24768905)$, $(-0.5225871,  -0.24768905)$,\\ $(-0.5225871,  -0.24768905)$, $(-0.5225871,  -0.24768905)$. 

The end point here is in the domain $x+y<0$, but is not close to any critical point of $f$. We find that the value of $f(-0.5225871,  -0.24768905)-M_0$ $=$ $-0.06196899168272818$. Therefore, we can choose $M_1=M_0-0.061968992$ $=$ $-0.061968992$ as the next lower bound estimate for $\min f$ in the concerned domain.  

Continue with $M_1=$ $-0.06196899168272818$: The random initial point is $(0.67180663, -1.47719597)$. Running Backtracking New Q-Newton's method for $G_{M_1}$, we get the sequence of points:\\ $(-0.21324179, -0.40600826)$, $(-0.44163373, -0.22503746)$, $(-0.60779332, -0.29120854)$,\\ $(-0.67567597, -0.30618429)$, $(-0.68489123, -0.30625665)$, $(-0.68507292, -0.30623803)$,\\ $(-0.68507296, -0.30623802)$, $(-0.685073,   -0.30623802)$, $(-0.685073   ,-0.30623802)$.  

Again,  the end point here is in the domain $x+y<0$, but is not close to any critical point of $f$. We find that the value of $f(-0.685073   ,-0.30623802)-M_1$ $=$ $-0.01060545418333153$. Therefore, we can choose $M_2=M_1-0.01060545418333153$ $=$ $-0.07257444619$ as the next lower bound estimate for $\min f$ in the concerned domain.    

Continue with $M_2=-0.07257444619$: The random initial point is $(0.68565678, -0.90941219)$. Running Backtracking New Q-Newton's method for $G_{M_2}$, we get the sequence of points:\\ $(-0.1515453,  -0.07688046)$, $(-0.23533302, -0.11778859)$, $(-0.36993082, -0.18298619)$,\\ $(-0.5502725,  -0.26540423)$, $(-0.67244411, -0.30841356)$, $(-0.70462466, -0.31280385)$,\\ $(-0.70678697, -0.31271293)$, $(-0.70679752, -0.31271166)$. 
   
Again,  the end point here is in the domain $x+y<0$. Now, it gets fairly close to the global minimizer of $f$ in this domain. We find that the value of $f(-0.70679752, -0.31271166)-M_2$ $=$ $-0.00015342244591928789$. Therefore, we can choose $M_3=M_2-0.00015342244592$ $=$ $-0.07272786863$ as the next lower bound estimate for $\min f$ in the concerned domain.  

  Continue with $M_3=-0.07272786863$: The random initial point is $(-0.65260407 -0.16301343)$. Running Backtracking New Q-Newton's method for $G_{M_2}$, we get the sequence of points:\\ $(-0.66547817, -0.28644322)$, $(-0.70420061, -0.31143069)$, $(-0.70709159, -0.31279723)$,\\ $(-0.70710672, -0.31280114)$. 
  
  Again, the end point is in the domain $x+y<0$, and indeed it is very close to the global minimizer in the domain. We observe that the value of $f(-0.70710672, -0.31280114)-M_3$ is, while still $<0$, very small. Hence, we can conclude (correctly) that the end point is a global minimizer. We check this by running one more time with the new value $M_4=M_3-3.1e-08$. 
  
  Running Backtracking New Q-Newton's method for the function $G_{M_4}$, from the randomly chosen initial point $(0.37597523, -0.91234048)$, we get the sequence of points:\\ $(-0.44234965, -0.06819544)$, $(-0.5447063,  -0.19358355)$, $(-0.66585017, -0.28453739)$,\\ $(-0.70420049, -0.31125032)$, $(-0.7070916,  -0.31279619)$, $(-0.70710678, -0.31280116)$. 
  
  \subsection{Example 6} In this example, we test our methods with a linear programming problem. The problem is as follows: Find $arg\min_{(x,y)\in S}f(x,y)$ where $f(x,y)=-40x-30y$, and $$S=\{(x,y)\in \mathbf{R}^2:~x+y\leq 12,~2x+y\leq 16, ~x\geq 0, ~y\geq 0 \}.$$
  
  The global minimum is $-400$, obtained for $(x,y)=(4,8)$. 
  
  Linear programming are usually difficult for iterative methods (in particular, Newton's type method), since the gradient of the cost function is a constant and the Hessian of the cost function is $0$. 
  
  For this problem, we find that the method of dividing with $d((x,y),\partial S)^2$ can lead to sequences which cross $\partial S$. We find that with the method of changing the function $f$ to be a big constant on $\mathbf{R}^2\backslash S$ indeed can produce sequences which approximate the global minimizers. Below we report experiments for the new function $g(x,y)=f(x,y)$ if $(x,y)\in S$, and $g(x,y)=1000$ if $(x,y)\in \mathbf{R}^2\backslash S$. 

  Running Backtracking New Q Newton's method for the function $g(x,y)$, with initial point $(0.1,0.1)$ we obtain the sequence of points (all in S): $(4.56428571, 3.44821429)$,  $(5.68035714, 4.28526786)$, $(5.75011161, 4.33758371)$, $(5.75011161, 4.33758371)$, $(5.80072885, 4.39576044)$, $(5.80159205, 4.39583457)$, $(5.80202368, 4.3958531)$, $(5.80205066, 4.39585426)$, $(5.80206414, 4.39585484)$, $(5.80207089, 4.39585513)$, $(5.80207173, 4.39585516)$, $(5.80207215, 4.39585518)$, $(5.80207236, 4.39585519)$, $(5.80207239, 4.39585519)$. The function value at the last point is $f(5.80207239, 4.39585519)=-363.95855134613464
$. 
  
   Running Backtracking New Q Newton's method for the function $g(x,y)$, with initial point $(1,2)$ we obtain the sequence of points (all in S): $(3.23214286, 3.67410714)$, $(4.34821429, 4.51116071)$, $(4.34832006, 6.18830041)$, $(4.37053239, 6.88502747)$, $(4.4401263,  6.95855197)$, $(4.48607869, 6.99847415)$, $(4.48773112, 7.02416628)$, $(4.48774402 7,.02436737)$, $(4.48775047, 7.02446792)$, $(4.48775209, 7.02449305)$, $(4.48775219, 7.02449462)$, $(4.48775224, 7.02449541)$, $(4.48775224, 7.02449551)$. The function value at the last point is $f(4.48775224, 7.02449551)=-390.2449550156782
$. 
   
 Running Backtracking New Q Newton's method for the function $g(x,y)$, with initial point $(2,1)$ we obtain the sequence of points (all in S):  $(4.23214286, 2.67410714)$, $(3.82457748, 7.44440135)$, $(4.01653468, 7.71777055)$, $(4.07413609, 7.77070227)$, $(4.09326204, 7.78987838)$, $(4.09733609, 7.79377355)$, $(4.09771509, 7.80059851)$, $(4.09873397, 7.80157261)$, $(4.09878135, 7.80242612)$, $(4.09878172, 7.80243279)$, $(4.09878191, 7.80243612)$, $(4.09878191, 7.80243618)$. The function value at the last point is $f(4.09878191, 7.80243618)=-398.0243616122982$, extremely close to the global minimum value!
 
 \subsection{Example 7} In this example we consider a constrained problem with a quadratic cost function \cite{gould-toint}. The question is: Find $arg\min _{(x,y)\in S}f(x,y)$ where $f(x,y)=-2(x-0.25)^2+2(y-0.5)^2$, and $$S=\{(x,y)\in \mathbf{R}^2:~x+y\leq 1, ~6x+2y\leq 3,~x,y\geq 0\}.$$
 
 This has a global minimum at the point $(0,0.5)$, with function value $f(0,0.5)=-0.125$. 
 
  In \cite{farrell-etal}, the tunnelling/deflation method was used to solve the associated equations coming from the Karush-Kuhn-Tucker optimality conditions. We can also use our method to solve these equations, similar to above. Here, we illustrate how to use  our method directly for the cost function $f$, without going through  the Karush-Kuhn-Tucker optimality conditions.
  
  We consider a new cost function $F(x,y)=f(x,y)$ if $(x,y)\in S$, and $F(x,y)=1000$ when $(x,y)\notin S$. Here are some runs. 
  
  Running Backtracking New Q-Newton's method for the function $F(x,y)$, with initial point $(0.1,0.1)$, we get a sequence:  $(0.025, 0.35)$, $(0.00625, 0.35)$, $(0.00117188, 0.35)$, $(0.00104637, 0.5)$, $(2.68394854e-04, 5.00000000e-01)$. The last point is very close to the global minimum in $S$. 
  
    Running Backtracking New Q-Newton's method for the function $F(x,y)$, with initial point $(0.2,0.3)$, we get a sequence:  $(0.15, 0.5)$, $(0.05, 0.5)$, $(0.025, 0.5)$, $(0.0109375, 0.5)$, $(0.0034668, 0.5)$, $(0.00154076, 0.5)$, $(0.00076432, 0.5)$, $(0.00076432, 0.5)$. Again, the last point is very close to the global minimum. 
    
     Running Backtracking New Q-Newton's method for the function $F(x,y)$, with initial point $(0.24,0.48)$, we get a sequence: $(0.23, 0.5)$, $(0.21, 0.5)$, $(0.17, 0.5)$, $(0.09, 0.5)$, $(0.01, 0.5)$, $(0.00142857, 0.5 )$, $(4.57589286e-04, 5.00000000e-01)$.  Again, the last point is very close to the global minimum.  
    
 \subsection{Example 8} In this example we find roots of a special function in a bounded domain. The question is: Find roots of the Bessel function $jv(1,z)$ inside the domain $S=\{z=x+iy:~-5\leq x,y\leq 5\}$. We use the library for special functions in Python to do computations with the Bessel function.  
 
 We consider the cost function $f(x,y)=|jv(1,x+iy)|^2/2$. We define a new cost function $F(x,y)=f(x,y)$ if $(x,y)\in S$, $F(x,y)=1000$ if $(x,y)\notin F$. Here are some experiments (initial points are randomly generated inside $S$). 
 
   Running Backtracking New Q-Newton's method for the function $F(x,y)$, with initial point $(3.61713097, 1.21693436)$, we get a sequence: $(3.81002318, 0.70240853)$, $(3.80130537, 0.23407724)$, $(3.82500859, 0.01478695)$, $(3.83166043e+00, 2.93087041e-05)$, $(3.83170597e+00, 3.48394778e-10)$, $(3.83170597e+00, 1.74197389e-10)$, $(3.83170597e+00, -6.20221574e-20)$. The last point is very close to a root of $jv(1,z)$ in $S$. 
 
  Running Backtracking New Q-Newton's method for the function $F(x,y)$, with initial point $(0.77926808, 3.75383432)$, we get a sequence: $(0.68962951, 3.20441492)$, $(0.22533733, 2.64615371)$, $(-0.28003512,  2.08563803)$, $(0.10219435, 1.52027522)$, $(-0.0071105,   0.95589959)$, $(7.82278286e-05, 4.21454021e-01)$, $(-4.48653604e-08,  6.18471202e-02)$, $(1.34987291e-14, 2.35499690e-04)$, $(-3.02313546e-13,  1.30608358e-11)$, $(-2.10926406e-12,  4.01821256e-13)$. The last point is very close to $(0,0)$, which is a root of $jv(1,z)$ inside $S$. 
  
  Running Backtracking New Q-Newton's method for the function $F(x,y)$, with initial point $(-2.1267499,  -0.96193073)$, we get a sequence: $(-2.44714686, -0.43462968)$, $(-3.44602951, -0.09106988)$, $(-3.64119905, -0.05036513)$, $(-3.82515578e+00, -2.52956231e-03)$, $(-3.83168871e+00, -4.33099427e-06)$, $(-3.83170597e+00, -1.95131964e-11)$, $(-3.83170597e+00, -9.75659827e-12)$. This is close to a root of $jv(1,z)$. 
  
  We have done many runs, and haven't seen a case where the sequence has a cluster point on $\partial S$. This is supported by Theorem \ref{Theorem}. In contrast, if we use the method of dividing $f$ by $d(.,\partial S)^2$, we often see the case where the constructed sequence converges to a root outside of $S$.  
  
 \subsection{Example 9} We revisit the function in Example 1 and Example 5.   Consider the following constrained optimization problem: 
 
$ arg\min _S f (x,y)$, where $S=\{(x,y)\in \mathbf{R}^2:~x+y\leq 0\}$ and $f(x,y)=-xye^{-x^2-y^2}+y^2/2$.

We consider a new cost function $F(x,y)=f(x,y)$ if $(x,y)\in S$, $F(x,y)=1000$ if $(x,y)\notin S$. 

Running Backtracking New Q-Newton's method for the function $F(x,y)$, with initial point $(0.5,  -0.5003)$, we get a sequence: $(-0.57530003, -0.28065)$, $(-0.64608507, -0.30173188)$,\\ $(-0.67661335, -0.30806708)$,$(-0.69177324, -0.31061912)$, $(-0.70704079, -0.31289215)$,\\ $(-0.70710676, -0.31280116)$. The last point is very close to the global minimum inside $S$. 

An interesting point here is that the point $(0.5,-0.5003)$ is in the basin of attraction of the global minimum {\bf outside of S} of the dynamics associated to Backtracking New Q-Newton's method applied to $f(x,y)$. By a simple change of the cost function, it belongs to the basin of attraction of the basin of attraction of the global minimum {\bf inside S} of the dynamics associated to Backtracking New Q-Newton's method applied to $F(x,y)$!

The way we run experiments here is much simpler that that in Example 5. There, we need to take care that the minimum of the function $f(x,y)-M$ is as close to $0$ as possible, and have to update the estimate $M$ for the global minimum of $f(x,y)$. 
   
If we apply the method Armijo's Backtracking line search with constraints, discussed in Section   2.2, to the cost function $f(x,y)$,  where the constraint condition is $x+y\leq 0$, to New Q-Newton's method, then we obtain a different sequence: $(0.5,-0.5003)$, $(0.49994643,-0.49999645)$, $(0.49993973, -0.49995854)$, $(0.49993638,-0.49993958)$, $(0.49993596,-0.49993721)$,\\ $(0.49993575,-0.49993603)$, $(0.4999357,-0.49993573)$, $(0.49993569,-0.4993569)$. The function value at the last point is approximately $0.276581$.  Here, instead of converging to the global minimum inside the domain $x+y< 0$, the constructed sequence -while being inside the constrained domain - converges to a point on the boundary $x+y=0$.   
   
 \section{Conclusions} In this paper we introduced two new methods to avoid a given closed set $A$. The first method is to divide the cost function $f$ by $d(x,A)^N$ for a suitable exponent $A$. This method is suitable in case one wants to avoid known points, to hop to a new component of the solution set, or to solve constrained optimization problems with no (local) minima on $A$. The second method is to change the value of $f$ on $A$ to a big constant. This method is more suitable for constrained optimization with (local) minima on the boundary. Experiments illustrate that the new methods are promising. Combination of these methods with each other, or with other methods can yield better performances. For example, our method can be used to quickly find a good initial point for Linear Programming problems, at which Linear Programming methods can be used to find exact global minima.  
   
{\bf A strategy for constrained optimization}: From the experiments above, we observe that the methods proposed in the current paper works well with constrained optimization problems for which the global minima lie strictly inside the interior of the constrained domain. In the case where the global minima lie on the boundary of the domain, here we present an idea on how to proceed. We know that when applying Backtracking New Q-Newton's method to the cost function $F(x)$ (which equals $f(x)$ when constraints are satisfied, and equals a big constant $M$ when constraints are not satisfied), then any cluster point of the constructed sequence - when belongs to the interior of the constrained domain - must be a critical point of $f$ and not a saddle point of $f$. Hence, we check at these cluster points first. In addition, if for many random choices of the initial point, we found that the constructed sequence tends to the boundary, then it may be an indication that we should explore the restriction of the function on the boundary. 

In case the constraints include both inequalities and equalities, e.g. $h_1(x)=\ldots =h_j(x)=0$ and $h_{j+1}(x),\ldots, h_l(x)\leq 0$,  one may need to first find a feasible point (or to check if the feasible set is empty). This can be done in 2 ways. 

Way 1:  We can add slack variables $y_{j+1},\ldots ,y_{l}$ and transform the constraints into the equivalent form $h_1(x)=\ldots =h_j(x)=0$ and $h_{j+1}(x)+y_{j+1}^2=\ldots =  h_l(x)+y_l^2= 0$. We can then use Backtracking New Q-Newton's method to find roots of this system of equations. 

Way 2: Assume that it is easy to find feasible points satisfying the conditions $h_{j+1}(x),\ldots, h_l(x)\leq 0$. We then apply Backtracking New Q-Newton's method to the cost function defined by $F(x)=h_1(x)^2+\ldots +h_j(x)^2$ if $h_{j+1}(x),\ldots, h_l(x)\leq 0$, and $F(x)=$ a big constant $M$ when the mentioned conditions are not satisfied.

\end{document}